\documentclass[10pt]{article}

\usepackage{cite}

\usepackage{amsfonts, amssymb, amsmath, amsthm, epsf}

\usepackage{graphicx}

\usepackage{fullpage}

\newtheorem{theorem}{Theorem}[section]
\newtheorem{lemma}{Lemma}[section]
\newtheorem{corollary}{Corollary}[section]
\newtheorem{condition}{Condition}[section]

\theoremstyle{definition}
\newtheorem{definition}{Definition}[section]
\newtheorem{example}{Example}[section]
\newtheorem{remark}{Remark}[section]

\numberwithin{equation}{section}
\numberwithin{figure}{section}

\newcommand{\const}{{\rm const}}

\renewcommand{\phi}{{\varphi}}

\newcommand{\cR}{{\mathcal R}}
\newcommand{\cH}{{\mathcal H}}

\newcommand{\oQ}{{\overline Q}}

\newcommand{\ob}{{\overline b}}
\newcommand{\ox}{{\overline x}}

\newcommand{\bbR}{{\mathbb R}}

\let\phi=\varphi
\newcommand{\al}{\alpha}
\newcommand{\ep}{\varepsilon}

\title{Reaction-diffusion equations with spatially distributed hysteresis}
\author{Pavel Gurevich
\footnote{Free Univerity of Berlin, Peoples' Friendship University
of Russia, email: gurevichp@gmail.com} , Roman Shamin
\footnote{Russian Academy of
Sciences, Far Eastern Branch, Institute of Marine Geology and Geophysics, email: roman@shamin.ru} , Sergey Tikhomirov
\footnote{Free Univerity of Berlin, email:
sergey.tikhomirov@gmail.com} }

\date{\today}

\begin{document}

\maketitle

\begin{abstract}
The paper deals with reaction-diffusion equations involving a hysteretic discontinuity in the source
term, which is defined at each spatial point. In particular, such problems describe chemical reactions and
biological processes in which diffusive and nondiffusive substances interact according to hysteresis law. We find sufficient conditions that guarantee the existence and uniqueness of solutions as well as their continuous dependence on initial data.
\end{abstract}






\section{Introduction}

\begin{figure}[ht]
\begin{minipage}[b]{0.53\linewidth}
        \centering
      \includegraphics[width=\linewidth]{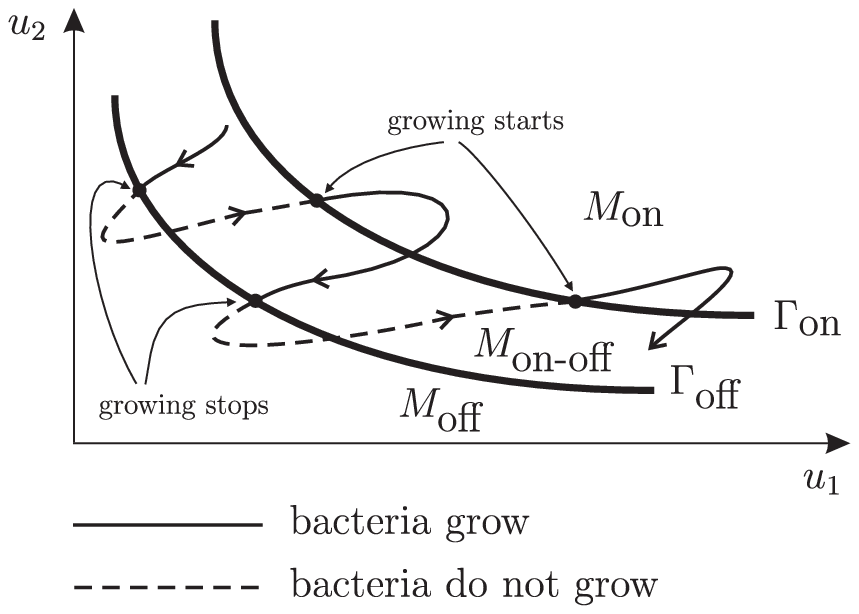}
        \caption{Regions of different behavior of hysteresis $\cH$}
        \label{fig1}
\end{minipage}
\hspace{0.5cm}
\begin{minipage}[b]{0.41\linewidth}
        \centering
        \includegraphics[width=\linewidth]{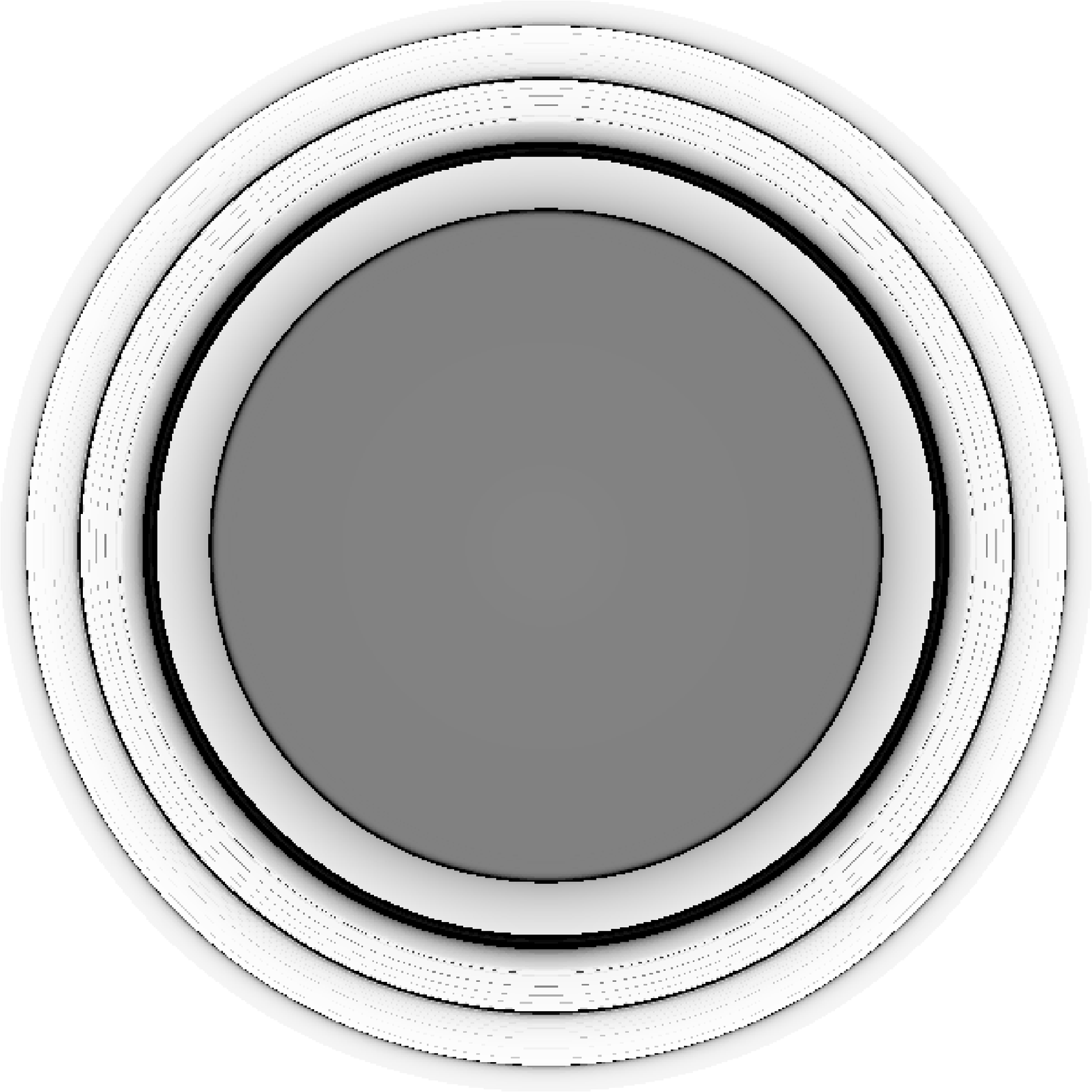}
        \caption{Density of bacteria after growth has stopped}
        \label{fig_BacteriaPattern}
        \end{minipage}
\end{figure}

The paper deals with reaction-diffusion equations involving a
hysteretic discontinuity which is defined at each spatial point.
In particular, these problems describe chemical reactions and
biological processes in which diffusive and nondiffusive
substances interact according to hysteresis law. We illustrate
this by a model describing a growth of a colony of bacteria
(Salmonella typhimurium) on a petri plate
(see.~\cite{Jaeger1,Jaeger2}). Let $Q\subset\bbR^n$ be a bounded
domain, $B(x,t)$ denote the density of nondiffusing bacteria in
$Q$, while $u_1(x,t)$ and $u_2(x,t)$ denote the concentrations of
diffusing buffer (pH level) and histidine (nutrient)  in $Q$,
respectively.  These three unknown functions satisfy the following
equations in $Q$:
\begin{equation}\label{eqBact1}
\left\{
\begin{aligned}
\dfrac{\partial B}{\partial t}&=avB,\\
\dfrac{\partial u_1}{\partial t}&=D_1\Delta u_1-a_1 vB,\\
\dfrac{\partial u_2}{\partial t}&=D_2\Delta u_2-a_2 vB\\
\end{aligned}\right.
\end{equation}
supplemented by initial and no-flux (Neumann) boundary conditions.
In those equations $D_1,D_2,a,a_1,a_2>0$ are given constants. The function
$v=v(x,t)$ corresponds to the growth rate of bacteria and is
defined by hysteresis law. In the simplest case, $v(x,t)$ takes
value $1$ if $u_1(x,t)$ and $u_2(x,t)$ are large enough and value
$0$ if $u_1(x,t)$ and $u_2(x,t)$ are small enough. More precisely,
one defines two curves $\Gamma_{\rm on}$ and $\Gamma_{\rm off}$ on
the plane $(u_1,u_2)$, which divide the first quadrant into three
parts $M_{\rm on}$, $M_{\rm off}$, and $M_{\text{on-off}}$. Now
$v(x,t)=1$ whenever $(u_1(x,t),u_2(x,t))\in M_{\rm on}$ and
$v(x,t)=0$ whenever $(u_1(x,t),u_2(x,t))\in M_{\rm off}$, while
$v(x,t)$ takes either value $1$ or $0$ in $M_{\text{on-off}}$
depending on whether the trajectory $(u_1(x,t),u_2(x,t))$ entered
the region $M_{\text{on-off}}$ through $\Gamma_{\rm on}$ or
$\Gamma_{\rm off}$. This hysteretic behavior is depicted at
Fig.~\ref{fig1}, see more details in~\cite{Jaeger1,Jaeger2}.

In the above example, hysteresis $v(x,t)$ may switch at different
spatial points at different time moments. This allows one to
divide the spatial domain $Q$ into connected subdomains: in each
of these subdomains, hysteresis $v(x,t)$ takes the same value ($1$
or $0$) and thus defines spatial topology of itself. The
boundaries between the subdomains are free boundaries whose motion
depends both on the reaction-diffusion equations and hysteresis.
The interplay between those two leads to formation of
spatio-temporal patterns. First numerical simulations exhibiting
such patterns have been carried out in~\cite{Jaeger1,Jaeger2} for
the case where $Q$ is a disc in $\mathbb R^2$. The appearing
pattern corresponds to concentric rings which are formed by
$B(x,t)$ as $t\to\infty$ (see Fig.~\ref{fig_BacteriaPattern}).

First rigorous results about the existence of solutions of
parabolic equations with hysteresis in the source term have been
obtained in~\cite{Alt, VisintinSpatHyst,Kopfova} for multi-valued
hysteresis. Formal asymptotic expansions of solutions were
recently obtained for some special case in~\cite{Ilin}. Questions
about the uniqueness of solutions and their continuous dependence
on initial data as well as a thorough analysis of pattern
formation remained open.

All the above questions are closely connected with slow-fast
system where hysteresis is replaced by a ``fast'' function
satisfying an ordinary differential equation with a small
parameter at the time derivative and (typically) cubic
nonlinearity. Then a possibility of singular perturbation limit is
generally an open question. Some results in this direction can be
found, e.g., in~\cite{Plotnikov, Evans}, where the authors study
equations of the form $u_t=\Delta\Phi(u)$ with cubic nonlinearity
$\Phi$. However, this direction is beyond the scope of our paper.

In this paper, we obtain theorems on   well-posedness of initial
boundary-value problems for reaction-diffusion equations with
discontinuous hysteresis in the source term. Our main assumption
is the so-called {\it spatial transversality} of initial data. In
the one-dimensional case, this roughly speaking means that the
initial function has nonvanishing derivative at the free boundary
between the above-mentioned subdomains. If this is the case, one
can guarantee that a solution exists on a small time interval and
can be continued as long as it remains spatially transverse. Moreover,
if solution is unique it continuously depends on the initial data.

Under further natural restrictions on hysteresis, it is possible
to prove that the solution is unique. For the completeness of
exposition, we formulate the corresponding theorem and refer the
reader to \cite{GurTikhUniq} for the proof.

For the clarity, all the results of the presented paper are proved for
a one-dimensional domain and a scalar reaction-diffusion equation.
But we note that the developed methods are not based on the
maximum principle; therefore, they can be
applied to systems of reaction-diffusion equations, including
those in multi-dimensional domains (see, e.g.,~\eqref{eqBact1}).

The paper is organized as follows. In Sec.~\ref{secSetting}, we
define functional spaces, introduce spatially distributed
hysteresis (i.e., defined at every spatial point) and set the
prototype problem:
\begin{equation}\label{eqMainProblem}
\left\{
\begin{aligned}
&u_t=u_{xx}+f(u,v),\quad  x\in(0,1),\ t>0,\\
&u_x|_{x=0}=u_x|_{x=1}=0,\\
&u|_{t=0}=\phi(x),
\end{aligned}\right.
\end{equation}
where $v(x,t)=\cH(u(x,\cdot))(t)$ represents hysteresis at the
point $x$. Loosely speaking, it is defined as follows. We fix two
thresholds $\alpha$ and $\beta$, $\alpha<\beta$ (analogues of
$\Gamma_{\rm on}$ and $\Gamma_{\rm off}$ in the above example).
Further, we fix two functions $H_1(u)$ and $H_2(u)$. Let
$x\in(0,1)$ be fixed. If $u(x,t)\le\alpha$, then we set
$v(x,t)=H_1(u(x,t))$; if $u(x,t)\ge\beta$, then we set
$v(x,t)=H_2(u(x,t))$; if $\alpha<u(x,t)<\beta$, then we set
$v(x,t)=H_1(u(x,t))$ or $H_2(u(x,t))$ depending on whether
$u(x,t)$ entered the interval $(\alpha,\beta)$ through $\alpha$ or
$\beta$, respectively (see Fig.~\ref{figHyst}).

Next, we discuss assumptions concerning the nonlinearity $f$, the
hysteresis operator $\cH$, and the transversality of the initial
function $\phi$. The most important notions here are the spatial
topology of hysteresis and transversality. We illustrate them by
the following situation. Suppose that, for a function $u(x,t)$,
there exists a continuous function $b(t)$, $t\in[0,T]$, taking
values in $(0,1)$ such that
\begin{equation}\label{eqIntroHystTopol}
\cH(u(x,\cdot))(t)=\begin{cases} H_1(u(x,t)), & 0\le x\le  b(t),\\
H_2(u(x,t)), & b(t)<x\le 1.
\end{cases}
\end{equation}
Then we say that $u$ {\it preserves spatial topology} of
hysteresis in the sense that, for all $t$, there are exactly two
subinterval $(0,b(t))$ and $(b(t),1)$; hysteresis is given by
$H_1(u)$ on one of them and by $H_2(u)$ on the other. In this
situation, $u$ is said to be {\it transverse} on $[0,T]$ if
\begin{enumerate}
\item $u(x,t)<\beta$ for $x\in[0,b(t)]$,

\item $u(x,t)>\alpha$ for $x\in(b(t),1]$

\item $u_x(b(t),t)>0$ whenever $u(b(t),t)=\alpha$.
\end{enumerate}

In the end of Sec.~\ref{secSetting}, we formulate the main results
of the paper: Theorem~\ref{thLocalExistence} (local existence of
transverse solutions preserving topology), Theorem~\ref{tCont}
(global existence of transverse solutions), and
Theorem~\ref{tContDepInitData} and
Corollary~\ref{corContDepInitData} (continuous dependence on
initial data). To make the exposure complete, we also formulate
Theorem~\ref{tUniqueness} (uniqueness of transverse solutions),
which is proved in~\cite{GurTikhUniq}.

We note that the failure of spatial topology may be caused by the
failure of transversality but need not necessarily be (see
Remark~\ref{remTransvNotPresSpTop}). In the latter case our
methods keep working.

In Sec.~\ref{SecAux}, we collect auxiliary results. First, we
recall a theorem on the well-posedness for linear parabolic
problems and then for semilinear problems. We show that the
hysteresis operator $\cH$ can be treated as a continuous operator
in suitable functional spaces, provided that its domain consists
of spatially transverse functions (although the function
$v(x,t)=\cH(u(x,\cdot))(t)$ still may have jumps at different
spatial points at different time moments). This allows us to prove
Theorem~\ref{tNP0}, which is the main tool in the study of
problem~\eqref{eqMainProblem}. This theorem deals with the
auxiliary problem
\begin{equation}\label{eqAuxProblem}
\left\{
\begin{aligned}
&u_t=u_{xx}+f(u,v_0),\quad  x\in(0,1),\ t>0,\\
&u_x|_{x=0}=u_x|_{x=1}=0,\\
&u|_{t=0}=\phi(x).
\end{aligned}\right.
\end{equation}
The meaning of $v_0$ is the following. Suppose we are given
functions $u_0(x,t)$ and $b_0(t)$. Then
\begin{equation}\label{eqIntroAuxHystTopol}
v_0(x,t)=\begin{cases} H_1(u_0(x,t)), & 0\le x\le  b_0(t),\\
H_2(u_0(x,t)), & b_0(t)<x\le 1.
\end{cases}
\end{equation}
One can think of $b_0(t)$ as of a free boundary defining spatial
topology of $v_0(x,t)$ in the sense similar to the above. However,
$v_0$ need not coincide with hysteresis $\cH(u_0)$. (It does only
if $u_0$ preserves spatial topology and the corresponding
subintervals are divided by the point $b_0(t)$.)

In particular, Theorem~\ref{tNP0} states that a solution $u$ of
problem~\eqref{eqAuxProblem} exists and is unique  on some time
interval $[0,T]$, where $T>0$ does not depend on $u_0$ and $b_0$.

In Sec.~\ref{SecProofs}, we prove the main results from
Sec.~\ref{secSetting}. The local existence is proved by means of
Theorem~\ref{tNP0} and the Schauder fixed-point theorem for a map
defined on pairs $(u_0,b_0)\in C([0,1]\times[0,T])\times C[0,T]$.
To prove global existence, we show that any local transverse
solution can be continued as long as it remains transverse. The
continuous dependence of solutions is also based on
Theorem~\ref{tNP0}.

In Sec.~\ref{secGeneral}, we generalize the dissipativity
condition for the nonlinearities $f(u,v)$ and $H_1(u),H_2(u)$ (cf.
Remark~\ref{remGeneral} and Example~\ref{Exh1h2}).

\section{Setting of the problem}\label{secSetting}

\subsection{Functional spaces}

Let $\Omega\subset\bbR^n$, $n\ge2$, be a domain with piece-wise
smooth boundary. We denote by $L_q(\Omega)$, $q>1$, the space with
the norm
$$
\|v\|_{L_q(\Omega)}=\left(\int_{\Omega} |v(y)|^q dy\right)^{1/q}.
$$

For an integer $l\ge 0$, we denote by $C^l(\overline\Omega)$ the
space of functions which are continuous and have continuous
derivatives up to order $l$ in $\overline\Omega$. For $l=0$, we
write $C(\overline\Omega)$. Denote  by
$C^\gamma(\overline\Omega)$, $0<\gamma<1$, the H\"older space with
the norm
$$
\|v\|_{C^\gamma(\overline\Omega)}=\|v\|_{C(\overline\Omega)}+\sup\limits_{x,y\in\Omega}\dfrac{|v(x)-v(y)|}{|x-y|^\gamma}.
$$

If $\Omega=(a,b)$, we will write $L_q(a,b)$, $C[a,b]$, $C^l[a,b]$,
and $C^\gamma[a,b]$. If $\Omega=(0,1)$, we will write $L_q$, $C$,
$C^l$, and $C^\gamma$, respectively. We also denote
$C^\infty=\bigcap\limits_{l=1}^\infty C^l$.

For a natural $l$, denote by $W_q^l=W_q^l(0,1)$ the space with the
norm
$$
\|v\|_{W_q^l}=\sum\limits_{j=0}^l \|v^{(j)}\|_{L_q},
$$
where $v^{(j)}$ is the generalized derivative of order $j$.

For any noninteger $l>0$, denote by $W_q^l=W_q^l(0,1)$ the space
with the norm
$$
\|v\|_{W_q^l}=\|v\|_{W_q^{[l]}}+\left(\int_0^1dx\int_0^1
\dfrac{|v^{([l])}(x)-v^{([l])}(y)|^q}{|x-y|^{1+q(l-[l])}}dy\right)^{1/q},
$$
where $[l]$ is the integer part of $l$.

Let $Q_T=(0,1)\times (0,T)$. Denote by $C^{1,0}(\overline Q_T)$
the space of function $u(x,t)$ such that $u,u_x\in C(\overline
Q_T)$.

We also introduce the anisotropic Sobolev space $ W_q^{2,1}(Q_T) $
with the norm
$$
\|u\|_{W_q^{2,1}(Q_T)} =\left(\int_0^T \|u(\cdot,t)\|_{W_q^2}^q\,
dt+ \int_0^T \|u_t(\cdot,t)\|_{L_q}^q\, dt\right)^{1/q}.
$$

The following embedding result holds (see Lemma 3.3 in~\cite[Chap.
2]{LadSolUral}).

\begin{lemma}\label{lEmbed}
If $u\in W_q^{2,1}(Q_T)$ with $q>3$, then, for any
$0\le\gamma<1-3/q$, $u,u_x\in C^\gamma(\overline Q_T)$ and
$$
\|u\|_{C^\gamma(\overline Q_T)}+\|u_x\|_{C^\gamma(\overline
Q_T)}\le c  \|u\|_{W_q^{2,1}(Q_T)},
$$
where $c>0$ depends on $q$, $\gamma$, and $T$ but does not depend
on $u$.
\end{lemma}
Throughout the paper, we fix $q$ and $\gamma$ such that
\begin{equation}\label{eqQGamma}
q>3, \quad 0<\gamma< 1-3/q.
\end{equation}
In what follows, we will consider solutions of parabolic problems
in the space $W_q^{2,1}(Q_T)$. To define a space of initial data,
we will use the fact that if $u\in W_q^{2,1}(Q_T)$, then the trace
$u|_{t=t_0}$ is well defined and belongs to $W_q^{2-2/q}$ for all
$t_0\in[0,T]$ (see Lemma 2.4 in~\cite[Chap. 2]{LadSolUral}).
Moreover, we can define the space $W_{q,N}^{2-2/q}$ as the
subspace of functions from $W_q^{2-2/q}$ with the zero Neumann
boundary conditions (it is well defined due to~\cite[Secs. 4.3.3
and 4.4.1 ]{Triebel}). Below we will also use the following embedding result

\begin{lemma}\label{lEmbed'}
If $\phi \in W_q^{2 - 2/q}$ with $q>3$, then, for any
$0\le\gamma<1-3/q$, $\phi, \phi'\in C^\gamma$ and
$$
\|\phi\|_{C^\gamma} +\|\phi'\|_{C^\gamma}\le c  \|u\|_{W_q^{2 - 2/q}},
$$
where $c>0$ depends on $q$ and $\gamma$ but does not depend
on $\phi$.
\end{lemma}

\subsection{Hysteresis}

In this section, we introduce a hysteresis operator defined for
functions of time variable $t$. Then we extend the definition to a
spatially distributed hysteresis acting on a space of functions of
time variable $t$ and space variable $x$.

We fix two numbers $\alpha$ and $\beta$ such that $\alpha<\beta$.
The numbers $\alpha$ and $\beta$ will play a role of thresholds
for the hysteresis operator. Next, we introduce continuous
functions
$$
H_1:(-\infty,\beta]\mapsto\bbR,\qquad H_2:[\alpha,\infty)
\mapsto\bbR.
$$
It is convenient to extend them to $\bbR$ as follows:
\begin{equation}\label{eqH12Extension}
\begin{aligned}
H_1(u)&=H_1(\beta)\ \text{for } u>\beta, \\
H_2(u)&=H_2(\alpha)\ \text{for } u<\alpha.
\end{aligned}
\end{equation}

We assume throughout that the following condition holds.

\begin{condition}\label{condH}
There is a number $\sigma\in(0,1]$ such that, for any $U>0$, there
exists $M=M(U)>0$ with the following property$:$
$$
|H_j(u)-H_j(\hat u)|\le M|u-\hat u|^\sigma,\quad j=1,2,
$$
whenever  $|u|,|\hat u|\le U$.
\end{condition}

We fix $T>0$ and denote by $C_r[0,T )$ the linear space of
functions which are continuous on the right in $[0,T )$. For any
$\zeta_0\in\{1,2\}$ ({\it initial configuration}) and $g\in
C[0,T]$ ({\it input}), we introduce the {\it configuration}
function
$$
\zeta:\{1,2\}\times C[0,T]\to C_r[0,T),\quad
\zeta(t)=\zeta(\zeta_0,g)(t)
$$
 as follows. Let
$X_t=\{t'\in(0,t]:g(t')=\alpha\ \text{or } \beta\}$. Then
$$
\zeta(0)=\begin{cases}
1 & \text{if } g(0 )\le \alpha,\\
2 & \text{if } g(0)\ge \beta,\\
\zeta_0 & \text{if } g(0)\in(\alpha,\beta)\\
\end{cases}
$$
and for $t\in(0,T]$
$$
\zeta(t)=\begin{cases}
\zeta(0) & \text{if } X_t=\varnothing,\\
1 & \text{if } X_t\ne\varnothing\ \text{and } g(\max X_t)=\alpha, \\
2 & \text{if } X_t\ne\varnothing\ \text{and } g(\max X_t)=\beta.
\end{cases}
$$

\begin{figure}[ht]
        \centering
        \includegraphics[height=120pt]{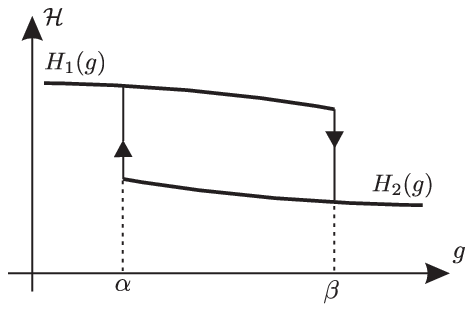}
        \caption{The hysteresis operator $\cH$}
        \label{figHyst}
\end{figure}

Now we introduce the {\it hysteresis operator\/}
(cf.~\cite{KrasnBook,Visintin})
$$
\cH: \{1,2\}\times C[0,T]\to C_r[0,T )
$$
by the following rule. For any initial configuration
$\zeta_0\in\{1,2\}$   and input $g\in C[0,T]$, the function
$\cH(\zeta_0,g):[0,T]\to\bbR$ ({\it output}) is given by
$$
\cH(\zeta_0,g)(t)=H_{\zeta(t)}(g(t)),
$$
where $\zeta(t)$ is the configuration function defined above (see
Fig.~\ref{figHyst}).


\begin{remark}
One usually assumes that $ H_1(u)- H_2(u)$ is sign-constant for
$u\in[\alpha,\beta]$. However, we will never use this assumption
in our paper.
\end{remark}

Now we introduce a {\it spatially distributed hysteresis}. Assume
that the initial configuration and the input function depend on
spatial variable $x\in[0, 1]$. Denote them by $\xi_0(x)$ and
$u(x,t)$, where
$$
\xi_0:[0, 1]\mapsto \{1, 2\},\qquad u:[0,1]\times[0,T]\mapsto\bbR.
$$

Let $u(x,\cdot)\in C[0,T]$. Denote $\phi(x)=u(x,0)$.

\begin{definition}\label{defGeneralConsistency}
We say that the  function $\phi(x)$ and the initial configuration
$\xi_0(x)$ are {\em consistent} if, for any $x\in[0,1]$,
$$ \xi_0(x)\in
\begin{cases}
\{1\} & \text{if } \phi(x)\le\alpha,\\
\{2\} & \text{if } \phi(x)\ge\beta,\\
\{1,2\} & \text{if } \phi(x)\in(\alpha,\beta).\\
\end{cases}
$$
\end{definition}

Assume that  $\xi_0(x)$  and   $\phi(x)=u(x,0)$ are consistent.
Then we can define the function
\begin{equation}\label{eqSDH}
v(x,t)=\cH(\xi_0(x),u(x,\cdot))(t),
\end{equation}
which is called {\it spatially distributed hysteresis}.

The {\it spatial configuration} is given by
\begin{equation}\label{eqSDC}
\xi(x,t)=\zeta(\xi_0(x),u(x,\cdot))(t).
\end{equation}
Note that the consistency of $\xi_0(x)$ and $\phi(x)$ and the fact
that $\xi(x,\cdot)$ is right-continuous guarantee that
$$
\lim\limits_{t\to 0}\xi(x,t)=\xi_0(x).
$$

In this paper, we shall deal will spatially distributed hysteresis
whose spatial configuration has finitely many discontinuity points
in $x$ for each $t$. In particular, we assume throughout that the
following holds.

\begin{condition}\label{condXiGeneral}
The initial configuration $\xi_0(x)$ has finitely many
discontinuity points in $(0,1)$, i.e., there are  points
$0=\overline b_0<\overline b_1<\dots<\overline b_M<\overline
b_{M+1}=1$ {\rm (}$M\ge 1${\rm)} such that
\begin{enumerate}
\item  $\xi_0(x)=\const $ for $x\in(\overline b_i,\overline
b_{i+1})$, $i=0,\dots,M$,

\item  $\xi_0(\overline b_i+0)\ne\xi_0(\overline b_i-0)$,
$i=1,\dots,M$.
\end{enumerate}
\end{condition}

Let us  introduce notions of transversality and   hysteresis
spatial
 topology, which will play a central role further on.

\begin{definition}\label{condTransverse}
We say that  a function $\phi \in C^1$ is {\em transverse} $($with
respect to a spatial  configuration $\xi_0(x))$ if it is
consistent with $\xi_0(x)$ and the following holds:
\begin{enumerate}
\item if $\phi(\ox)=\alpha$ and $\phi'(\ox)=0$ for some
$\ox\in[0,1]$, then $\xi_0(\ox)=1$ in a neighborhood of $\ox$;

\item if $\phi(\ox)=\beta$ and $\phi'(\ox)=0$ for some
$\ox\in[0,1]$, then $\xi_0(\ox)=2$ in a neighborhood of $\ox$.
\end{enumerate}
\end{definition}

\begin{definition}
We say that a function $u\in C^{1,0}(\overline Q_T)$ is {\em
transverse} {\em on} $[0,T]$ $($with respect to a spatial
configuration $\xi(x,t))$ if, for every fixed $t\in[0,T]$, the
function $u(\cdot,t)$ is transverse with respect to the spatial
configuration $\xi(\cdot,t)$.
\end{definition}

\begin{definition}
We say that a function $u\in C^{1,0}(\overline Q_T)$ {\em
preserves  spatial topology $($of a spatial configuration
$\xi(x,t))$ on $[0,T]$} if there is $M>0$ such that, for
$t\in[0,T]$, there are continuous functions
\begin{equation}\label{eqb0bM1}
0\equiv b_0(t)<b_1(t)<\dots< b_M(t)<b_{M+1}(t)\equiv 1
\end{equation}
 with the
properties
\begin{enumerate}
\item  $\xi(x,t)=\const$ for $x\in(b_i(t), b_{i+1}(t))$,
$i=0,\dots,M$,

\item  $\xi(b_i(t)+0,t)\ne\xi(b_i(t)-0,t)$, $i=1,\dots,M$.
\end{enumerate}
We will also say in this case that $u$ is {\em topology
preserving}.

\end{definition}

\begin{remark}\label{remTransvNotPresSpTop}
Neither of the above two definitions is equivalent to another. It
may happen that the number of discontinuity points $b_i(t)$ in the
interval $(0,1)$ decreases at some moment $t_1$, while $u(x,t)$
remains transverse with respect to $\xi(x,t)$ for each fixed
$t\in[0,T]$. This may happen if $b_i(t)-b_{i+1}(t)\to 0$ as $t\to
t_1$ for some $i\in\{0,\dots,M\}$ (see
Fig.~\ref{figChangeTopology}).

\begin{figure}[ht]
        \centering
        \includegraphics[width=\textwidth]{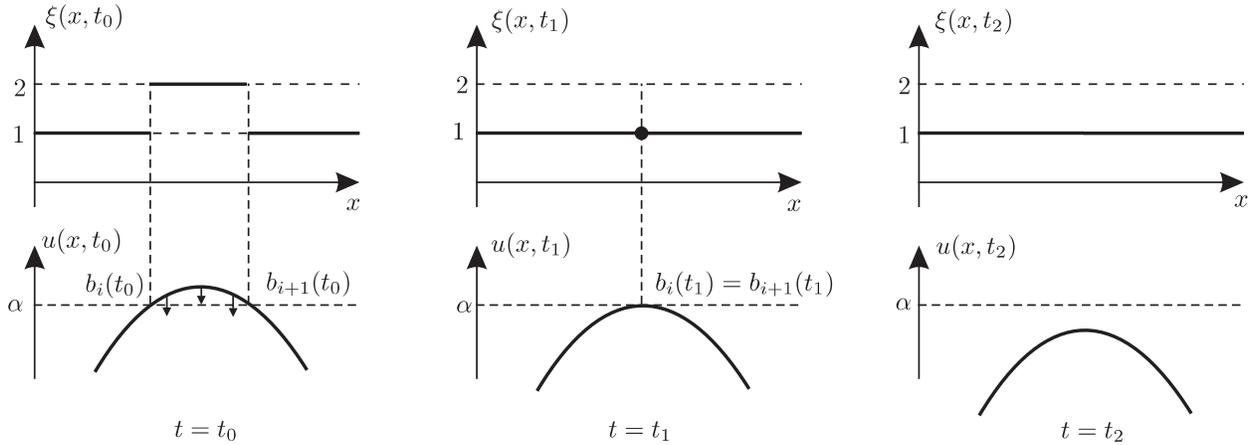}
        \caption{Topology changes while transversality condition holds}
        \label{figChangeTopology}
\end{figure}

On the other hand, the equalities $u(b_i(t),t)=\alpha$ or $\beta$
and $u_x(b_i(t),t)=0$ may hold for some $t$ and some
$i\in\{1,\dots,M\}$ without changing the spatial topology (see
Fig.~\ref{figNontransversal}).

\begin{figure}[ht]
        \centering
        \includegraphics[width=\textwidth]{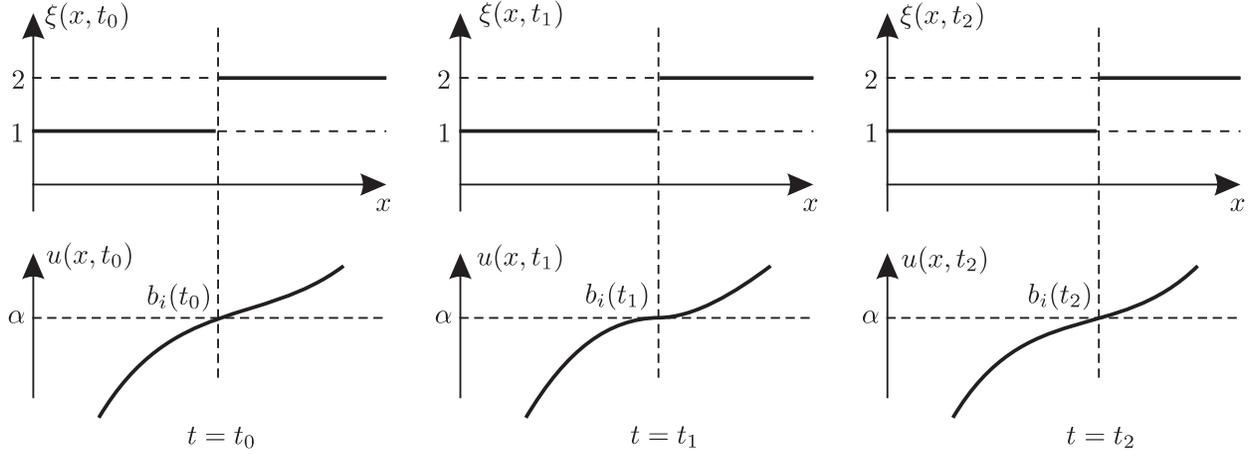}
        \caption{Topology preserves while transversality condition fails}
        \label{figNontransversal}
\end{figure}

\end{remark}

\subsection{Reaction-diffusion equations with hysteresis}

The main object of this paper is  the reaction-diffusion equation
\begin{equation}
u_t=u_{xx}+f(u,v),\qquad (x,t)\in Q_T,\label{eq1}
\end{equation}
where $v=v(x,t)$ represents the spatially distributed hysteresis
given by~\eqref{eqSDH}:
$$
v(x,t)=\cH(\xi_0(x),u(x,\cdot))(t).
$$
We also impose the Neumann boundary conditions
\begin{equation}
u_x|_{x=0}=u_x|_{x=1}=0\label{eq2}
\end{equation}
and the initial condition
\begin{equation}
u|_{t=0}=\phi(x),\quad x\in(0,1).\label{eq3}
\end{equation}

Along with Conditions~\ref{condH} and~\ref{condXiGeneral}, we
assume throughout that the right-hand side $f$, the hysteresis
branches $H_j(u)$, and the initial data satisfy the following
conditions.

\begin{condition}[Lipschitz continuity]\label{condRegularity}
For any bounded set $B\subset\mathbb R^2$, there is a constant
$L=L(B)>0$ such that
$$
|f(u_1,v_1)-f(u_2,v_2)|\le L(|u_1-u_2|+|v_1-v_2|)\quad \forall
(u_j,v_j)\in B,\ j=1,2.
$$
\end{condition}

\begin{condition}[dissipativity]\label{condDissip}
For any sufficiently large $U>0$ and for all  $u$  such that
$|u|\le U$, we have
$$
f(U,H_j(u))<0,\quad f(-U,H_j(u))>0,\quad  j=1,2.
$$
\end{condition}

\begin{remark}\label{remGeneral}
In Sec.~\ref{secGeneral}, we will generalize
Condition~\ref{condDissip} to the following:
\begin{enumerate}
\item for all sufficiently large $u$,
$$
f(u,H_2(u))\le 0,\quad f(-u,H_1(-u))\ge 0;
$$

\item there is a Lipschitz continuous function $h(u)$ such that
$uh(u)>0$ for $u\ne 0$ and, for any (small) $\mu>0$, there exists
$U_\mu >0$ such that the function
$$
f_\mu(u,v)=f(u,v)-\mu h(u)
$$
satisfies
$$
f_\mu(U_\mu,H_j(u))\le 0,\quad f_\mu(-U_\mu,H_j(u))\ge 0\quad
\forall |u|\le U_\mu,\ j=1,2.
$$
\end{enumerate}
\end{remark}

\begin{example}\label{Exh1h2}
The simplest example of $f(u,v)$ is
$$
f(u,v)=-d u+v,
$$
where $d>0$ and
$$
\lim\limits_{U\to+\infty}(-d U+\max\limits_{|u|\le
U}H_j(u))<0,\qquad \lim\limits_{U\to-\infty}(-d
U+\max\limits_{|u|\le U}H_j(u))>0.
$$

The simplest example of the hysteresis in this case is given by
$$
H_1(u)\equiv h_1,\quad H_2(u)\equiv h_2,
$$
where $h_1$ and $h_2$ are arbitrary real numbers without any sign
restrictions (see Fig.~\ref{figh1h2}).

Moreover, Remark~\ref{remGeneral} with $h(u)\equiv u$ shows that
if $h_1\le 0\le h_2$, then one can set $d=0$.
\begin{figure}[ht]
        \centering
        \includegraphics[width = \textwidth]{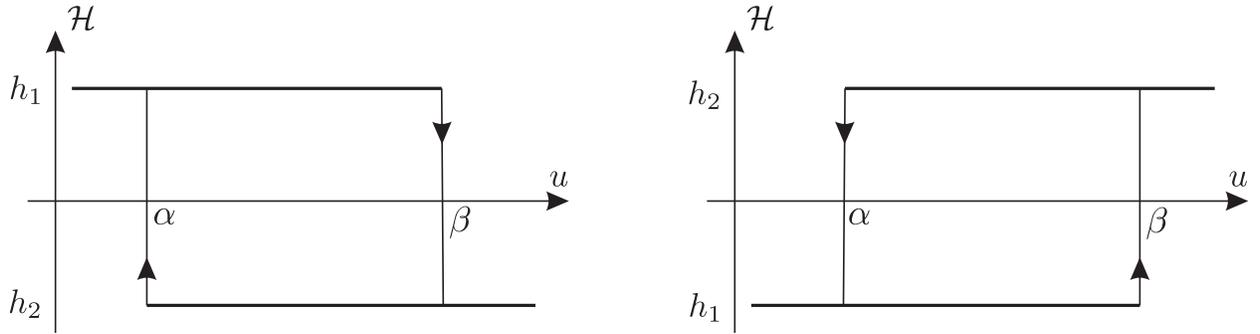}
        \caption{There are no sign restrictions for $h_1$ and $h_2$ in Example~\ref{Exh1h2}: both hysteresis operators are admissible,
        provided that $d>0$.}
        \label{figh1h2}
\end{figure}

\end{example}

The natural assumption on the initial data is as follows.

\begin{condition}[consistency]\label{condConsisten}
The initial spatial configuration $\xi_0(x)$  and the initial data
$\phi(x)$ are consistent in the sense of
Definition~$\ref{defGeneralConsistency}$.
\end{condition}

To prove the existence of a solution for
problem~\eqref{eq1}--\eqref{eq3}, we will   restrict a class of
admissible   initial data, namely we assume throughout the
following.

\begin{condition}\label{condTransversGeneral}
The initial function $\phi(x)$ is transverse in the sense of
Definition~$\ref{condTransverse}$.
\end{condition}

We give  a definition of a solution of
problem~\eqref{eq1}--\eqref{eq3},   assuming  that $\phi\in
W_{q,N}^{2-2/q}$.

\begin{definition}
A function $u(x,t)$ is called a {\em solution  of
problem~\eqref{eq1}--\eqref{eq3} {\rm (}in $Q_T${\rm) }} if $u\in
W_q^{2,1}(Q_T)$, $v(x,t)$ is measurable, $u$ and $v$ satisfy
equation~\eqref{eq1} for a.e. $(x,t)\in Q_T$, and
conditions~\eqref{eq2} and~\eqref{eq3} are satisfied in the sense
of traces.
\end{definition}

It follows from this definition and from Lemma~\ref{lEmbed} that
any solution $u(x,t)$ belongs to $C(\oQ_T)$. Therefore, the
function $v(x,t)$ is well defined by~\eqref{eqSDH} and belongs to
$L_\infty(Q_T)$.

\subsection{Main results}

Our main results are the following three theorems. In all of them,
we assume that
Conditions~$\ref{condH}$--$\ref{condTransversGeneral}$ hold  and
that $q$ and $\gamma$ satisfy~\eqref{eqQGamma}.

\begin{theorem}[local existence]\label{thLocalExistence}
There is a number $T>0$ such that the following holds.
\begin{enumerate}
\item Any solution $u\in W_q^{2,1}(Q_T)$ of
problem~\eqref{eq1}--\eqref{eq3}  in $Q_T$ is transverse and
preserves spatial topology. For each $t\in[0,T]$, the function
$v(\cdot,t)$  has exactly $M$ discontinuity  points
$b_1(t)<\dots<b_M(t)$ in $(0,1)$. Moreover,
\begin{enumerate}
\item $b_i(0)=\overline b_i$,

\item $b_i\in C^\gamma[0,T]$.
\end{enumerate}
\item There is at least one transverse topology preserving
solution $u\in W_q^{2,1}(Q_T)$ of problem~\eqref{eq1}--\eqref{eq3}
in~$Q_T$.
\end{enumerate}
\end{theorem}

The next theorem deals with continuation of solutions to their
maximal intervals of transverse existence.

\begin{definition} We say that $[0,T_{max})$, $T_{max}\le\infty$, is a {\em maximal
interval of transverse existence} of a solution $u$ of
problem~\eqref{eq1}--\eqref{eq3} if
\begin{enumerate}
\item for any $T<T_{max}$, the function $u$ is a transverse
solution of problem~\eqref{eq1}--\eqref{eq3}  in $Q_T$,

\item either $T_{max}=\infty$, or  $T_{max}<\infty$ and $u$ is a
solution in $Q_{T_{max}}$, but $u(\cdot, T_{max})$ is not
transverse with respect to $\xi(\cdot,T_{max})$.
\end{enumerate}
\end{definition}

Note that a  solution need not be topology preserving inside its
maximal interval of transverse existence  (see
Remark~\ref{remTransvNotPresSpTop}).

\begin{theorem}[continuation]\label{tCont}
Let $u\in W_q^{2,1}(Q_{t_0})$ be a transverse topology preserving
 solution of problem~\eqref{eq1}--\eqref{eq3}  in $Q_{t_0}$. Then
it can be continued to a maximal interval of transverse existence
$[0,T_{max})$, where $T_{max}>t_0$   may depend on continuation.
\end{theorem}

The following theorem shows that if a transverse topology
preserving solution is unique, then it continuously depends on
initial function $\phi$ and initial configuration~$\xi_0$.

\begin{theorem}[continuous dependence on initial data]\label{tContDepInitData}
We assume that the following hold.
\begin{enumerate}
\item There is a number $T>0$ such that
problem~\eqref{eq1}--\eqref{eq3} with initial function $\phi\in
W_{q,N}^{2-2/q}$ and initial configuration $\xi_0(x)$ defined by
its discontinuity points $\ob_1<\dots<\ob_M$ admits a unique
transverse topology preserving solution $u\in W_q^{2,1}(Q_s)$ in
$Q_s$ for any $s\le T$.

\item Let $\phi_n\in W_{q,N}^{2-2/q}$, $n=1,2,\dots$, be a
sequence of other initial functions such that
$\|\phi-\phi_n\|_{W_{q,N}^{2-2/q}}\to 0$ as $n\to\infty$.

\item Let $\xi_{0n}(x)$,  $n=1,2,\dots$, be a sequence of other
initial configurations defined by their discontinuity points
$\ob_{1n}<\dots<\ob_{Mn}$ such that $\xi_{0n}(x)=\xi_0(x)$ for
$x\in (0,\min(b_1,b_{1n}))$ and $\ob_{jn}-\ob_j\to 0$ as
$n\to\infty$, $j=1,\dots,M$.
\end{enumerate}
Then, for all sufficiently large $n$,
problem~\eqref{eq1}--\eqref{eq3} with initial function $\phi_n$
and initial configuration $\xi_{0n}$ has at least one transverse
topology preserving solution $u_n\in W_q^{2,1}(Q_T)$. Each
sequence of such solutions satisfies
$$
\|u_n-u\|_{W_q^{2,1}(Q_T)}\to 0,\quad \|b_{jn}-b_j\|_{C[0,T]}\to
0,\quad j=1,\dots,M,\quad n\to\infty,
$$
where $b_j(t)$ and $b_{jn}(t)$ are the respective discontinuity
points of the configuration functions $\xi(t)$ and $\xi_j(t)$,
$j=1,\dots,M$.
\end{theorem}

\begin{figure}[ht]
        \centering
        \includegraphics[width=\textwidth]{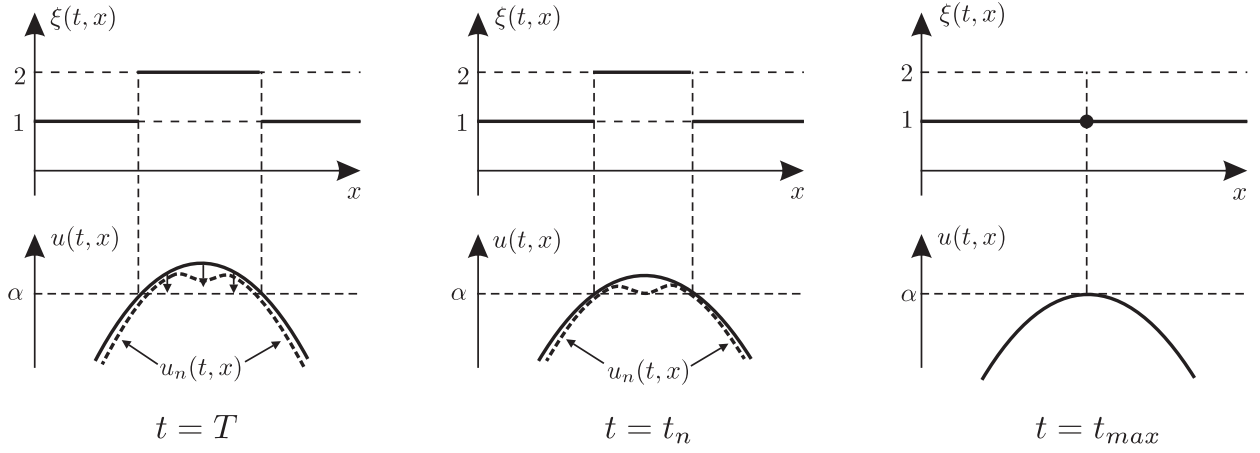}
        \caption{Approximation $u_n$ becomes nontransverse at a moment $t_n\in(T,t_{max})$}
        \label{fig5}
\end{figure}

As we have mentioned before, a solution $u$ need not be topology
preserving on its maximal interval of transverse existence
$[0,T_{max})$. But Theorem~\ref{tContDepInitData} deals only with
topology preserving solutions. The reason is the following. Let
$t_{max}$ ($t_{max}<T_{max}$)  be a moment at which the solution
$u$ changes topology. Although Theorem~\ref{tContDepInitData}
guarantees that $u_n$ remain transverse and topology preserving on
any interval $[0,T]\subset[0,t_{max})$ and approximate $u$, it may
happen that each $u_n$ becomes nontransverse at some moment
$t_n\in(T,t_{max}]$ for all $n\ge N=N(T)$. This situation is
illustrated by Fig.~\ref{fig5}.

However, if we a priori know that all $u_n$ are transverse on some
interval $[0,T]\subset[0,T_{max})$ (possibly with $T\ge t_{max}$),
then we can prove that $u_n$ approximate $u$. Now the topology of
$u$ can change. This happens if some neighboring discontinuity
points $b_j(t)$ and $b_{j+1}(t)$ of $\xi(x,t)$ converge to the
same point $\hat b_j$ as $t\to t_{max}$. Then, for $t>t_{max}$,
the point $b_j(t)=b_{j+1}(t)=\hat b_j$ is not a discontinuity
point of $\xi(x,t)$ any more. However, if we consider the
functions $b_j(t)$ and $b_{j+1}(t)$ on the whole interval $[0,T]$
assuming them constant for $t>t_{max}$, then we can prove that
$b_{jn}(t)$ approximate $b_j(t)$ in $C[0,T]$. Here we use the same
convention about $b_{jn}(t)$.

 Let us formulate this assertion as a
corollary from Theorem~\ref{tContDepInitData}.

\begin{corollary}\label{corContDepInitData}
\begin{enumerate}
\item We fix some time interval $[0,T]$, initial function $\phi\in
W_{q,N}^{2-2/q}$, and initial configuration $\xi_0(x)$ defined by
its discontinuity points $\ob_1<\dots<\ob_M$. Assume that
problem~\eqref{eq1}--\eqref{eq3} with initial data $\phi$ and
$\xi_0(x)$ admits a unique transverse solution $u\in
W_q^{2,1}(Q_s)$ in $Q_s$ for any $s\le T$.

\item Let assumptions $2$ and $3$ of
Theorem~$\ref{tContDepInitData}$ hold. Moreover, suppose that, for
each $n$, there exists a transverse solution $u_n\in
W_q^{2,1}(Q_T)$ of problem~\eqref{eq1}--\eqref{eq3} with initial
data $\phi_n$ and $\xi_{0n}(x)$.
\end{enumerate}
Then $$ \|u_n-u\|_{W_q^{2,1}(Q_T)}\to 0,\quad
\|b_{jn}-b_j\|_{C[0,T]}\to 0,\quad j=1,\dots,M,\quad n\to\infty.
$$
\end{corollary}

For the completeness of exposition, we also formulate the
uniqueness theorem, (see \cite{GurTikhUniq} for the proof).

In the uniqueness theorem, along with Conditions
$\ref{condH}$--$\ref{condTransversGeneral}$ we assume that the
following holds.
\begin{condition}\label{condUniq}
There is a number $\sigma\in[0,1)$ such that, for any $U>0$, there
exists $M=M(U)>0$ with the properties
\begin{equation}\label{eqCondH1}
|H_1(u)-H_1(\hat u)|\le \dfrac{M}{(\beta-u)^\sigma+(\beta-\hat
u)^\sigma}|u-\hat u|,\quad \forall u,\hat u\in[-U,\beta),
\end{equation}
\begin{equation}\label{eqCondH2}
|H_2(u)-H_2(\hat u)|\le \dfrac{M}{(u-\alpha)^\sigma+(\hat
u-\alpha)^\sigma}|u-\hat u|,\quad \forall u,\hat u\in(\alpha,U].
\end{equation}
\end{condition}

\begin{theorem}[uniqueness]\label{tUniqueness}
Let the functions $H_j(u)$, $j = 1, 2$, additionally satisfy
Condition~$\ref{condUniq}$. Assume that $u,\hat u\in
W_q^{2,1}(Q_{T_0})$ are two transverse solutions of
problem~\eqref{eq1}--\eqref{eq3} in $Q_{T_0}$ for some $T_0$. Then
$u=\hat u$.
\end{theorem}

\begin{remark}
\begin{enumerate}
  \item Any locally Lipschitz continuous functions $H_1(u)$ and
  $H_2(u)$ satisfy Condition~\ref{condUniq}. Moreover, this condition
  covers the important case where non-Lipschitz hysteresis branches $H_1(u)$
  and $H_2(u)$ appear as a singular-perturbation limit of slow-fast system (see \cite{GurTikhUniq}).

  \item On the other hand, any $H_1(u)$ and $H_2(u)$ satisfying
  Condition~\ref{condH} are locally H\"older continuous with
  exponent $1-\sigma$ on $(-\infty,\beta]$ and $[\alpha,\infty)$,
  respectively. Therefore, they satisfy Condition \ref{condH}.
\end{enumerate}
\end{remark}


\subsection{Technical simplification}\label{subsecSimple}

\begin{figure}[ht]
\begin{center}

\begin{minipage}[t]{0.30\linewidth}
        \centering
        \includegraphics[width=\linewidth]{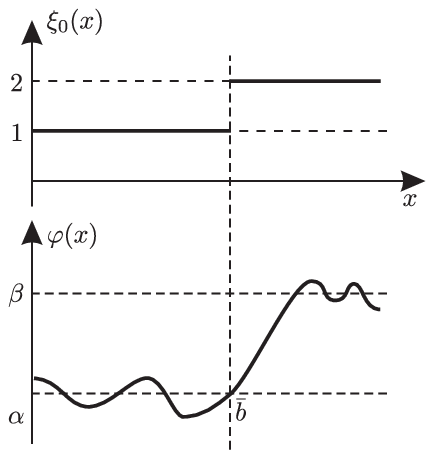}
        \caption{Initial data  satisfying Condition~\ref{condA1'}}
        \label{fig11}
\end{minipage}
\hspace{0.5cm}
\begin{minipage}[t]{0.45\linewidth}
        \centering
        \includegraphics[width=\linewidth]{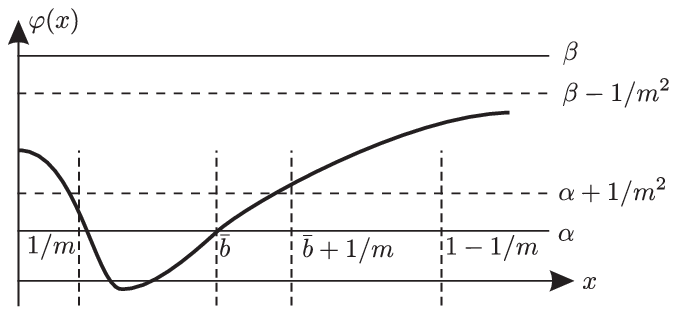}
        \caption{Initial data belongs to $E_m$}
        \label{fig6}
        \end{minipage}
        \end{center}
\end{figure}

For the clarity of exposition, we give detailed proofs of the main
results for the functions $\phi(x)$ and $\xi_0(x)$ satisfying the
following condition  (see Fig. \ref{fig11}).

Fix some $\ob\in (0,1)$.

\begin{condition}\label{condA1'}
\begin{enumerate}
\item For  $\ob\in(0,1)$, one has
\begin{equation}\label{eqxi0}
\xi_0(x)=\begin{cases}
1, & x\le\ob,\\
2, & x>\ob.
\end{cases}
\end{equation}
\item $\phi(x) < \beta$ for $x \in [0,\ob]$.

\item $\phi(x) > \alpha$ for $x \in (\ob,1]$.

\item If $\phi(\ob)=\alpha$, then $\phi'(\ob)>0$.

\end{enumerate}
\end{condition}


It follows from this condition that the hysteresis~\eqref{eqSDH}
at the initial moment is given by
$$
v|_{t=0}=\begin{cases}
H_1(\phi(x)), & x\le\ob,\\
H_2(\phi(x)), & x>\ob.
\end{cases}
$$

Clearly, Condition~\ref{condXiGeneral} (with $M=1$ and $\overline
b_1=\ob$) and Conditions~\ref{condConsisten}
and~\ref{condTransversGeneral} are satisfied in this case.

We denote by $E_m$, $m\in\mathbb N$, the set of pairs
$(\phi,\xi_0)$ satisfying Condition \ref{condConsisten} such that $\phi\in W_{q,N}^{2-2/q}$, $\xi_0(x)$ is
of the form~\eqref{eqxi0} and the following hold (see
Fig.~\ref{fig6}):
\begin{enumerate}
\item $\ob \in [1/m, 1 - 1/m]$,

\item $\phi(x)\le\beta-1/m^2$ for $x \in [0, \ob]$,

\item $\phi(x) \ge \alpha + 1/m^2$ for $x  \in [\ob + 1/m, 1]$,

\item if $x\in [\ob, \ob + 1/m]$ and $\phi(x) \in [\alpha, \alpha
+ 1/m^2]$, then $\phi'(x)\ge 1/m$,
%

\item $\|\phi\|_{W_{q,N}^{2-2/q}}\le m$.
\end{enumerate}


Note that item 4 in the definition of $E_m$ yields
\begin{equation}\label{Add13.1}
\varphi(x) \geq \al + \frac{1}{m}(x- \bar{b}), \quad x \in [\bar{b}, \bar{b} + 1/m].
\end{equation}
This observation easily implies $E_m\subset E_{m+1}$. Moreover, Lemma~\ref{lEm} below
shows that the unity of all sets $E_m$ coincides with the set of
all transverse data satisfying Condition~$\ref{condA1'}$. These sets allow us to measure the ``level of transversality'' of
the data. The precise statement is as follows.
\begin{lemma}\label{lEm}
\begin{enumerate}
\item Functions $\phi(x)$ and $\xi_0(x)$ satisfy
Condition~$\ref{condA1'}$ if and only if
$(\phi,\xi_0)\in\bigcup\limits_{m=1}^\infty E_m$.

\item Let $\phi_m\in W_{q,N}^{2-2/q}$, and let $\xi_m(x)$ be
defined analogously to~\eqref{eqxi0} with some $\ob_m$ instead of
$\ob$. If
\begin{enumerate}
\item $(\phi_m,\xi_m)\in E_m\setminus E_{m-1}$, $m=2,3,\dots$,

\item $\|\phi_m-\phi\|_{W_{q,N}^{2-2/q}}\to 0$ as $m\to \infty$
for some $\phi\in W_{q,N}^{2-2/q}$,

\item $ \ob_m-\ob\to0$ as $m\to \infty$ for some $\ob\in[0,1]$,
\end{enumerate}
then $\ob\in\{0,1\}$ or $\phi(x)$ is not transverse with respect
to $\xi_0(x)$, where $\xi_0(x)$ is given by~\eqref{eqxi0}.

\item Let $(\phi,\xi_0)\in E_m$ for some $m\in \mathbb N$. Then
there is $\varepsilon=\varepsilon(m)>0$ such that
$(\tilde\phi,\tilde\xi_0)\in E_{m+1}$ whenever
$\|\tilde\phi-\phi\|_{W_{q,N}^{2-2/q}}\le\varepsilon$, $|\tilde
b-\ob|\le\varepsilon$, and $\tilde\xi_0$ is given analogously
to~\eqref{eqxi0} with $\tilde b$ instead of $\ob$.
\end{enumerate}
\end{lemma}
\proof
1.  Let $\phi(x)$ and $\xi_0(x)$ satisfy
Condition~\ref{condA1'}.  Items 1, 2, 5 in the definition of $E_m$
directly follow from Condition~\ref{condA1'} for a
sufficiently large $m$.

Further, if $\phi(\ob)\ne\alpha$, then one can choose $m$ such
that   $\phi(x) \notin [\alpha, \alpha + 1/m^2]$ for all points
$x\in [\ob, \ob + 1/m]$. In this case item 4 in the definition of
$E_m$ becomes void and item 3 in the definition of $E_m$ is a trivial consequence of item 3 in Condition \ref{condA1'}.

Finally, if $\phi(\ob)=\alpha$, then $\phi'(\ob) > 0$ and hence
there exists $m>0$ such that $\phi'(x) \geq 1/m$ for $x \in [\ob,
\ob +1/m]$. This inequality implies item 4 of the definition of
$E_m$ and inequalities (\ref{Add13.1}). Applying item 3 of
Condition \ref{condConsisten}, we can further increase $m$ so that
item 3 hold.

If $(\phi,\xi_0)\in E_m$ for some $m$, then it is obvious that Condition \ref{condA1'} holds.

\medskip

2. Assertion 2 will follow from assertions 1 and 3. Indeed, if
$\ob\in(0,1)$ and $\phi(x)$ is transverse with respect to
$\xi_0(x)$, then assertion 1 implies that  $(\phi,\xi_0)\in
E_{m_0}$ for a sufficiently large $m_0$. Therefore, by
assertion~3, $(\phi_m,\xi_m)\in E_{m_0+1}$ for all sufficiently
large $m$, which contradicts   assumption (a) on the sequence
$(\phi_m,\xi_m)$.

\medskip

3. { \rm Let us prove assertion 3.}
Lemma \ref{lEmbed'} implies that
\begin{equation}\label{1.1}
\|\tilde{\phi} - \phi\|_{C} \leq c\ep, \quad \| \tilde{\varphi}' - \phi' \|_{C} \leq c\ep, \quad \|\tilde{\phi}'\|_{C^{\gamma}} \leq c(m+1),
\end{equation}
where $c>0$ does not depend on $\ep$ and $m$. Items 1, 2, 5 of the definition of $E_{m+1}$ are direct consequence of (\ref{1.1}), provided that $\ep>0$ is small enough.

Let us prove item 4 in the definition of $E_{m+1}$. Consider $x
\in [\tilde{b}, \tilde{b} +1/(m+1)]$ such that $\tilde{\phi}(x)
\in [\al, \al + 1/(m+1)^2]$. If $x \geq \ob$, then item 4 is a
direct consequence of \eqref{1.1}. Assume $x < \ob$. Then $|x -
\ob| \leq \ep$, and inequalities (\ref{1.1}) imply
$$
\tilde{\phi}(\ob) \leq \al + \frac{1}{(m+1)^2} + |x-\ob|c(m+1) \leq \al + \frac{1}{(m+1)^2} + \ep c(m+1).
$$
For small enough $\ep$, this implies that $\phi(\ob) \leq \al +
1/m^2$ and hence $\phi'(\ob) \geq 1/m$. This implies
$\tilde{\phi}'(\ob) \geq 1/m - \ep$ and (by inequalities
(\ref{1.1}))
$$
\tilde{\phi}'(x) \geq \frac{1}{m} - \ep - c(m+1)\ep^{\gamma}.
$$
The latter inequality implies $\tilde{\phi}'(x) \ge 1/(m+1)$ for
small enough $\ep$, which completes the proof of item 4.

Let us  now prove item 3 in the definition of $E_{m+1}$. Consider
$x \in [\tilde{b} + 1/(m+1), 1]$. If $x \geq \ob + 1/m$, then the
inequality $\phi(x) \geq \al + 1/m^2$ and \eqref{1.1} imply
$\tilde{\phi}(x)\geq 1/(m+1)^2$ for small enough $\ep$. Let $x <
\ob +1/m$. The inequality $|\tilde{b} - \ob| \leq \ep$ implies $x
- \ob \geq 1/(m+1) - \ep$.

 Therefore,
$$
\phi(x) \ge \al +\frac{1}{m}\left( \frac{1}{m+1} - \ep \right).
$$
Now, using \eqref{1.1}, we conclude that $\tilde{\phi}(x) \geq \al
+ 1/(m+1)^2$, provided that $\ep > 0$ is small enough, which
completes the proof of assertion 3.
\endproof

We will also need the following auxiliary statement.

\begin{lemma}\label{lEmA}
Let $(\phi, \xi_0) \in E_m$ and $\ob \in (0, 1)$ be such that
\eqref{eqxi0} holds. Then, for any functions $\psi_1, \psi_2 \in
C^1$  satisfying
\begin{equation}\label{14*}
\|\psi_i - \phi\|_{C} + \|\psi_i' - \phi'\|_C \leq
\dfrac{1}{4m^2}, \quad i = 1, 2,
\end{equation}
the following holds.
\begin{enumerate}
\item The equation $\psi_i(x) = \al$ $(i=1,2)$ has no more than
one root on the interval $[\ob, 1]$. \item If such a root exists,
let us denote it by $a_i;$ otherwise, we set $a_i = \ob$. Then
    \begin{equation}\label{14.1}
    a_1, a_2 \in [\ob, \ob + 1/m],
    \end{equation}
    \begin{equation}\label{3.1}
    |a_1 - a_2| \leq 2m \|\psi_1 - \psi_2\|_C.
    \end{equation}
\end{enumerate}
\end{lemma}
\begin{proof}
Item 1 and inclusions \eqref{14.1} are a trivial consequences of items 3 and 4 in the definition of $E_m$.

Let us prove item 2. First, let us assume that both $a_1, a_2$
satisfy $\psi_{1, 2}(a_{1, 2}) = \al$ (see Fig. \ref{figa1a2}).
Without loss of generality, we can assume that $a_2 > a_1$. Let us
fix an arbitrarily $x \in [a_1, a_2]$. Item 4 in the definition of
$E_m$ and \eqref{14*} imply that $\psi_2(x) \leq 0$ and $\psi_1(x)
\leq \al + 1/(2m^2)$. Hence, due to \eqref{14*}, $\phi(x) \le \al
+ 1/m^2$, which implies $\phi'(x) \geq 1/m$ and, thus, $\psi'_1(x)
\geq 1/(2m)$. The latter inequality yields
$$
|a_2 - a_1| \leq 2m|\psi_1(a_2) - \psi_2(a_2)|,
$$
which proves (\ref{3.1}).
\begin{figure}[ht]
\begin{center}
\begin{minipage}[t]{0.45\linewidth}
        \centering
        \includegraphics[width=\linewidth]{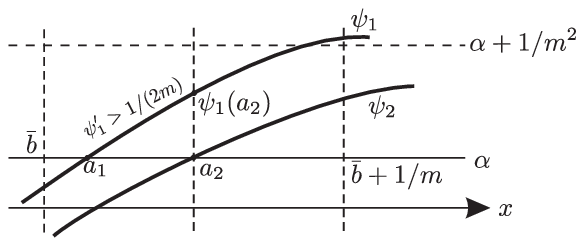}
        \caption{Both roots of $\psi_{1, 2}(x) = \al$ exist}
        \label{figa1a2}
\end{minipage}
\hspace{0.5cm}
\begin{minipage}[t]{0.45\linewidth}
        \centering
        \includegraphics[width=\linewidth]{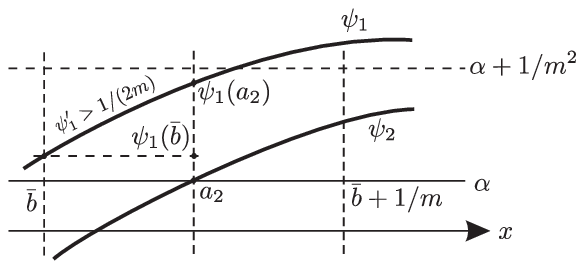}
        \caption{A root of $\psi_2(x) = \al$ exists; $a_1 = \ob$}
        \label{figa1a2'}
        \end{minipage}
        \end{center}
\end{figure}

Assume that $a_2$ satisfies $\psi_2(a_2) = \al$ and the equation $\psi_1(x) = \al$ has no roots on the interval $[\ob, 1]$ (see Fig. \ref{figa1a2'}).
 Using inequalities (\ref{Add13.1}) and \eqref{14*}, we conclude that $\psi_1(\ob +1/(2m))> 0$ and hence $\psi_1(\ob) > 0$. Arguing similarly to the previous case, we conclude that $\psi'_1(x) \geq 1/(2m)$ for $x \in [\ob, a_2]$ and hence
 $$
 |a_2 - a_1| = |a_2 - \ob| \leq 2m |\psi_1(a_2) - \psi_1(\ob)| \leq 2m |\psi_1(a_2)- \psi_2(a_2)|,
 $$
 which proves (\ref{3.1}).

If neither of the equations $\psi_{1, 2}(x) = \al$ has a root on
the interval $[\ob, 1]$, then $a_1 = a_2 = \ob$ and (\ref{3.1}) is
trivial.
\end{proof}

\section{Auxiliary results}\label{SecAux}

\subsection{Linear parabolic problem}

In this subsection, we formulate a well-known result on the
solvability of a linear parabolic equation in the anisotropic
Sobolev space $W_q^{2,1}(Q_T)$.

We consider the initial boundary-value problem
\begin{equation}\label{eqF}
\left\{
\begin{aligned}
& u_t=u_{xx}+F(x,t),\quad (x,t)\in Q_T,\\
& u_x|_{x=0}=u_x|_{x=1}=0,\\
& u|_{t=0}=\phi(x),\quad x\in(0,1).
\end{aligned}
\right.
\end{equation}
Combining the results of~\cite[Chap. 4]{LadSolUral} (including
Lemma~\ref{lEmbed} formulated above in our paper) with Theorem 3.1
in~\cite{Sobolevskii} and with the interpolation theory
in~\cite[Secs. 1.14.5, 4.3.3, and 4.4.1]{Triebel}, we obtain the
following result.

\begin{theorem}\label{tLinearHeat}
\begin{enumerate}
\item Assume that $q$ and $\gamma$ satisfy~\eqref{eqQGamma}. Fix
 numbers $T_0\ge T>0$. Let $F\in L_q(Q_T)$ and $\phi\in
W_{q,N}^{2-2/q}$. Then problem~$\eqref{eqF}$ has a unique solution
$u\in W_q^{2,1}(Q_T)$ and
$$
\|u\|_{W_q^{2,1}(Q_T)}+\max_{t\in[0,T]}\|u(\cdot,t)\|_{W_{q,N}^{2-2/q}}\le
c_1(\|\phi\|_{W_{q,N}^{2-2/q}}+\|F\|_{L_q(Q_T)}),
$$
$$
\|u\|_{C^\gamma(\overline Q_T)}+\|u_x\|_{C^\gamma(\overline
Q_T)}\le c_2 (\|\phi\|_{W_{q,N}^{2-2/q}}+\|F\|_{L_q(Q_T)}),
$$
The constants $c_1,c_2>0$ depend on $q$, $\gamma$, and $T_0$, but
do not depend on $T$, $\phi$, and $F$.
\end{enumerate}
\end{theorem}

\subsection{Semilinear parabolic problem}\label{subsecMild}

In this subsection, we consider an auxiliary semilinear initial
boundary-value problem
\begin{equation}\label{eqPf0}
\left\{
\begin{aligned}
& u_t=u_{xx}+f_0(u,x,t),\quad (x,t)\in Q_T,\\
& u_x|_{x=0}=u_x|_{x=1}=0,\\
& u|_{t=0}=\phi(x),\quad x\in(0,1),
\end{aligned}
\right.
\end{equation}

\subsubsection{$E_{\infty,T}$-mild solutions}

First, we  assume that $\phi\in L_\infty$. We also assume that the
function $f_0$ satisfies the following.
\begin{enumerate}
\item $f_0(u,x,t)$ is measurable with respect to
$(x,t)\in[0,1]\times[0,\infty)$ for all $u\in\mathbb R$.

\item For every bounded set $B\subset \mathbb R\times
[0,1]\times[0,\infty)$, there exists a constant $L=L(B)>0$ such
that
$$
|f_0(u,x,t)|\le L\quad \forall (u,x,t)\in B,
$$
$$
|f_0(u,x,t)-f_0(v,x,t)|\le L|u-v|\quad \forall (u,x,t),(v,x,t)\in
B.
$$
\end{enumerate}

We give one result from~\cite{Rothe} on the so-called
$E_{\infty,T}$-mild   solutions of  problem~\eqref{eqPf0}.

First, we   define the operator $A_0: D(A_0)\subset L_p\to L_p$
for $p>1$ by
$$
D(A_0)=\{\psi\in C^2: \psi'(0)=\psi'(1)=0\},
$$
$$
A_0\psi=-\psi_{xx}+\psi\quad \forall\psi\in D(A_0).
$$
The operator $A_0$ has a closure $A_p$ in $L_p$. The operators
$A_p$ are generators of analytic semigroups $S_p(t)$ in $L_p$. By
Lemma~1 in~\cite[p.~15]{Rothe}, $S_{p_1}(t)\subset S_{p_2}(t)$ for
any $p_1,p_2\in(1,\infty)$, $p_1\ge p_2$. By Lemma~2
in~\cite[p.~19]{Rothe},
$$
\sup\limits_{t\in[0,T]}\|S_p(t)\psi\|_{L_\infty}\le
\|\psi\|_{L_\infty}\quad \forall \psi\in L_\infty.
$$
Therefore, one can define the operators $P(t)$ as the restrictions
of $S_p(t)$ to the space $L_\infty$. These operators do not depend
on $p$ and are continuous from $L_\infty$ to $L_\infty$.
Furthermore, since
$$
P(t_1+t_2)=P(t_1)P(t_2)\quad \forall t_1,t_2\in[0,\infty),
$$
they define a semigroup in $L_\infty$.

Note that it is not a strongly continuous semigroup in $L_\infty$
because, due to Lemma~2 in~\cite[p.~19]{Rothe}, $P(t)\psi\to\psi$
in $L_\infty$ as $t\to 0$ if and only if $\psi\in C$.

\begin{definition}
Let $T\in (0,\infty]$. A {\em $E_{\infty,T}$-mild solution} of
problem~\eqref{eqPf0} for initial data $\phi\in L_\infty$ on the
time interval $[0,T)$ is a measurable function $u(x,t)$,
$(x,t)\in(0,1)\times(0,T)$, satisfying
$$
u(\cdot,t)\in L_\infty,\qquad \sup\limits_{s\in
(0,t)}\|u(\cdot,s)\|_{L_\infty}<\infty\quad \forall t\in(0,T),
$$
$$
u(\cdot,t)=P(t)\phi+\int_0^t
P(t-s)(f_0(u(\cdot,s),\cdot,s)+u(\cdot,s))\,ds\quad \forall
t\in(0,T),
$$
where the integral is an absolutely converging Bochner inegral in
$L_\infty$.
\end{definition}

\begin{definition}
We say that $T\in(0,\infty)$ is a {\em maximal existence time} for
a given initial data $\phi$ if problem~\eqref{eqPf0} has  a
$E_{\infty,T}$-mild solution on the interval $[0,T)$, but for any
$T'>T$, it has no $E_{\infty,T'}$-mild solution on the interval
$[0,T')$.
\end{definition}

The following lemma is formulated as Theorem 1 in~\cite[p.
111]{Rothe}.

\begin{lemma}\label{lRothe}
Assume that $\phi\in L_\infty$ and the above conditions on $f_0$
are satisfied. Then there exists a maximal existence time
$T\in(0,\infty]$ and  problem~\eqref{eqPf0} has a unique
$E_{\infty,T}$-mild solution on the interval $[0,T)$.

If the maximal existence time $T$ is finite, then
$$
\lim\limits_{t\to T}\|u(\cdot,t)\|_{L_\infty}=\infty.
$$
\end{lemma}

\subsubsection{Uniformly bounded solutions}

Now we formulate a result which states that if
problem~\eqref{eqPf0} has a solution $u\in W_q^{2,1}(Q_T)$, then
it is bounded uniformly with respect to $T\in(0,\infty)$. It can
be proved by regularizing the right-hand side and applying the
invariant-rectangles method, which can also be exploit for systems
of reaction-diffusion equations (see, e.g.,~\cite{Smoller}).

\begin{lemma}\label{lIR}
Let $\phi\in W_{q,N}^{2-2/q}$, $u\in W_q^{2,1}(Q_T)$, and
$F(x,t)=f_0(u(x,t),x,t)$ belongs to $L_q(Q_T)$. We assume that $u$
satisfies~\eqref{eqPf0} and the following hold for some $U>0$.
\begin{enumerate}
\item $f_0(\cdot,x,t)$ is continuous at the points $\pm U$
uniformly with respect to $(x,t)\in \oQ_T$ $($perhaps, after
modification on a subset of $\oQ_T$ of measure zero$)$,

\item $f_0(U,x,t)<0$, $f_0(-U,x,t)>0$ for a.e. $(x,t)\in \oQ_T$,

\item $\|\phi\|_C<U$.
\end{enumerate}
Then $\|u\|_{C(\oQ_T)}<U$.
\end{lemma}

\subsection{Semilinear parabolic problem with a special
nonlinearity}

In this subsection, we specify the nonlinearity $f_0(u,x,t)$
in~\eqref{eqPf0}.

Let us fix  initial data $\phi$ and $\xi_0$ satisfying
Condition~\ref{condA1'}. By Lemma~\ref{lEm}, $(\phi,\xi_0)\in E_m$
for a sufficiently large $m$. Let us fix such an $m$ (see
Fig~\ref{fig6}).
 Next, we choose a number $U>0$ such that
Condition~\ref{condDissip} holds for $U$ and
\begin{equation}\label{eqphi1}
\|\phi\|_C <U.
\end{equation}

Further, we fix a number $T>0$ and consider functions $u_0\in
C(\oQ_T)$ and $b_0\in C[0,T]$ with the following properties:
\begin{equation}\label{equ01}
\|u_0\|_{C(\oQ_T)}\le U,
\end{equation}
\begin{equation}\label{equ02}
b_0(t)\in [\ob,\ob+1/m], \quad t\in[0,T].
\end{equation}

Now we define the function $f_0(u,x,t)$ by
\begin{equation}\label{eqf0}
f_0(u,x,t)=f(u,v_0(x,t)),
\end{equation}
\begin{equation}\label{eqv0}
v_0(x,t)=\begin{cases} H_1(u_0(x,t)), & 0\le x\le  b_0(t),\\
H_2(u_0(x,t)), & b_0(t)<x\le 1,
\end{cases}
\end{equation}
where we use convention~\eqref{eqH12Extension}.

Note that, in general, $v_0(x,t)\ne
\cH(\xi_0(x),u_0(x,\cdot))(t)$. However, the following is true and
will be essentially used later on.
\begin{remark}\label{rH=v0}
Suppose that the following hold for each $t\in[0,T]$.
\begin{enumerate}
\item $u_0(x,0)=\phi(x)$.

\item The equation $ u(x,t) = \alpha$, on the interval $[\ob,1]$
has no more than one root. Let us define a function $a_0(t)$ as
follows: if the above root exists, then $a_0(t)$ is equal to this
root, and $a_0(t)=\ob$ otherwise.

\item The function
\begin{equation}\label{eqb0maxa0}
b_0(t)=\max\limits_{s\in[0,t]}a_0(s)
\end{equation}
satisfies $b_0(t)\in[\ob,\ob+1/m]$.

\item The equation $u_0(x,t)=\beta$ on the interval $[0,b_0(t)]$
has no roots.
\end{enumerate}
Then the nonlinearity $v_0(x,t)$  given by~\eqref{eqv0}  coincides
with the spatially distributed hysteresis defined for $u_0(x,t)$:
$$
v_0(x,t)\equiv \cH(\xi_0(x),u_0(x,\cdot))(t).
$$
\end{remark}

The next lemma establishes basic properties of the
``$\max$''-operator in~\eqref{eqb0maxa0}.

\begin{lemma}\label{lab}
Let $\lambda=0$ or $\gamma$. Then the following hold.
\begin{enumerate}
\item If $a\in C^\lambda[0,T]$ and $b(t)=\max\limits_{s\in[0,t]}
a(s)$, then $b\in C^\lambda[0,T]$ and
$$
\|b\|_{C^\lambda[0,T]}\le\|a\|_{C^\lambda[0,T]}.
$$

\item If $a_j\in C[0,T]$ and $b_j(t)=\max\limits_{s\in[0,t]}
a_j(s)$, $j=1,2$, then
$$
\|b_1-b_2\|_{C[0,T]}\le \|a_1-a_2\|_{C[0,T]}.
$$
\end{enumerate}
\end{lemma}
\proof We leave details to the reader.
\endproof

In what follows, we will need the continuous dependence of $v_0$
defined by~\eqref{eqv0} on $u_0$ and $b_0$. Denote by $\cR$ the
set of pairs $(u_0,b_0)\in C(\oQ_{T})\times C[0,T]$ satisfying
conditions~\eqref{equ01} and~\eqref{equ02}.

\begin{lemma}\label{lu0b0v0}
For any $(u_0,b_0),(\hat u_0,\hat b_0)\in \cR$ $(m\in\mathbb N)$,
let $v_0$ be defined by~\eqref{eqv0} and $\hat v_0$
by~\eqref{eqv0}, where $u_0$ and $b_0$ are replaced by $\hat u_0$
and $\hat b_0$, respectively. Then, for any $p\in[1,\infty)$,
\begin{equation}\label{equ0b0v0}
\|v_0-\hat v_0\|_{L_p(Q_T)}\le c_0 \left(T^{1/p} \|u_0-\hat
u_0\|_{C(\oQ_T)}^\sigma+\|b_0-\hat b_0\|_{L_1(0,T)}^{1/p}\right),
\end{equation}
where $\sigma$ is the constant in Condition~$\ref{condH}$ and
$c_0>0$ depends on $U$ and $p$, but does not depend on
$u_0,b_0,T$.
\end{lemma}
\proof We fix some $t\in[0,T]$ and assume that $b_0(t)\le \hat
b_0(t)$ for this $t$. Then, using~\eqref{eqv0} and omitting the
arguments of the integrands, we have
$$
\begin{aligned}
\int\limits_0^1 |v_0-\hat v_0|^p\,
dx&=\int\limits_0^{b_0(t)}|H_1(u_0)-H_1(\hat
u_0)|^p\,dx+\int\limits_{\hat b_0(t)}^1|H_2(u_0)-H_2(\hat
u_0)|^p\,dx\\
& +\int\limits_{b_0(t)}^{\hat b_0(t)}|H_2(u_0)-H_1(\hat
u_0)|^p\,dx.
\end{aligned}
$$
Using Condition~\ref{condH} in the first two integrals and the
boundedness of $H_j(u)$ for $|u|\le U$ in the second integral, we
obtain
$$
\int\limits_0^1 |v_0-\hat v_0|^p\, dx\le k_1(\|u_0-\hat
u_0\|_{C(\oQ_T)}^{\sigma p}+k_2|b_0(t)-\hat b_0(t)|),
$$
where $k_1>0$ depends on $U$ and $p$, but does not depend on
$u_0,b_0,T$.

Integrating the latter inequality with respect to $t$ from $0$ to
$T$ yields~\eqref{equ0b0v0}.
\endproof

Denote
\begin{equation}\label{eqVfU}
 V=\max\limits_{|u|\le
U}\{H_1(u),H_2(u)\},\qquad f_U=\max\limits_{|u|\le U, |v|\le V}
|f(u,v)|.
\end{equation}

The main step in finding a solution of
problem~\eqref{eq1}--\eqref{eq3}  is the following theorem.

\begin{theorem}\label{tNP0}
Let $f_0$ be defined by~\eqref{eqf0}, and let $q$ and $\gamma$
satisfy~\eqref{eqQGamma}. We fix an arbitrary $T_0>0$. Then the
following hold.
\begin{enumerate}
\item Problem~\eqref{eqPf0} has a unique solution $u\in
W_q^{2,1}(Q_T)$ and, for any $T\le T_0$,
\begin{equation}\label{eqNP02**}
\|u\|_{C(\oQ_T)}< U,
\end{equation}
\begin{equation}\label{eqNP02*}
\begin{aligned}
\|u\|_{W_q^{2,1}(Q_T)}+\max_{t\in[0,T]}\|u(\cdot,t)\|_{W_{q,N}^{2-2/q}} \le c_1(\|\phi\|_{W_{q,N}^{2-2/q}}+f_U),\\
\|u\|_{C^\gamma(\oQ_T)}+\|u_x\|_{C^\gamma(\oQ_T)}  \le c_2
(\|\phi\|_{W_{q,N}^{2-2/q}}+f_U),
\end{aligned}
\end{equation}
where $f_U$ is given by~\eqref{eqVfU} and $c_1,c_2>0$ depend only
on $T_0$ and do not depend on $m,u_0,b_0,\phi,u,T$.

\item The solution of problem~\eqref{eqPf0} continuously depends
on $\phi$, $u_0$ and $b_0$. In other words, if $u_n\in
W_q^{2,1}(Q_T)$, $n=1,2,\dots$, are   solutions of
problem~\eqref{eqPf0} with $\phi,u_0,b_0,v_0$ replaced by $\phi_n,
u_{0n},b_{0n},v_{0n}$, then
\begin{equation}\label{equ-un}
\|u_n-u\|_{W_q^{2,1}(Q_T)}\to 0,\quad n\to\infty,
\end{equation}
whenever
\begin{equation}\label{equ0-u0n}
\|\phi_n-\phi\|_{W_{q,N}^{2-2/q}}+\|u_{0n}-u_0\|_{C(\oQ_T)}+\|b_{0n}-b_0\|_{C[0,T]}\to
0,\quad n\to\infty.
\end{equation}

\item There is a number $T=T(m)\in (0,T_0]$ and a natural number
$N=N(m,U)\ge m$ which do not depend on $u_0,b_0,\phi,u$, such
that, for any $t\in[0,T]$, the following is true.
\begin{enumerate}
\item The equation $u(x,t)=\alpha$ on the interval $[\ob,1]$ has
no more than one root. If this root exists, we denote it by
$a(t);$ otherwise, we set $a(t)=\ob$ $($similarly to Lemma
$\ref{lEmA})$. In this case, $a(t)\in[\ob,\ob+1/N]$, $a\in
C^\gamma[0,T]$, and
\begin{equation}\label{equ2'}
\|a\|_{C^\gamma[0,T]}\le a^*,
\end{equation}
where $a^*>0$ depends on $m$, but does not depend on
$u_0,b_0,\phi$.

\item The hysteresis $\cH(\xi_0,u)$  and its configuration
function $\xi(x,t)$ have  exactly one discontinuity point $b(t);$
moreover, $b(t)=\max\limits_{s\in[0,t]}a(s)$ and
 $b\in C^\gamma[0,T]$.

\item  $(u(\cdot,t),\xi(\cdot,t))\in E_N$.
\end{enumerate}
\end{enumerate}
 \end{theorem}
\proof 1. Throughout the proof, we assume that $u_0$ and $b_0$ are
extended as continuous functions to $[0,T_0]$ in such a way
that~\eqref{equ01} and~\eqref{equ02} hold on this interval.

It follows from the definition~\eqref{eqf0} of the function $f_0$
and from the Lipschitz continuity of the function $f$
(Condition~\ref{condRegularity}) that the function $f_0(u,x,t)$
satisfies     assumptions 1 and 2 in Sec.~\ref{subsecMild}.
Therefore, by Lemma~\ref{lRothe}, there is $T_1\in(0, T_0]$ such
that problem~\eqref{eqPf0} has a unique $E_{\infty,T_1}$-mild
solution. Hence, $f_0(u(x,t),x,t)$ is in $L_\infty(Q_{T_1})$, and
Theorem~\ref{tLinearHeat} yields $u\in W_q^{2,1}(Q_{T_1})$.

Now we claim that
\begin{equation}\label{eqNP0uU}
\|u\|_{C(\oQ_{T_1})}<U.
\end{equation}
Indeed, the function $v_0(x,t)$ in~\eqref{eqf0} is bounded for
$(x,t)\in\oQ_{T_1}$; therefore, $f_0(u,x,t)=f(u,v_0(x,t))$
satisfies assumption~1 in Lemma~\ref{lIR} due to the uniform
continuity of $f(u,v)$ on bounded sets. Assumption~2 in
Lemma~\ref{lIR} holds due to Condition~\ref{condDissip} and the
choice of $u_0$ in~\eqref{equ01}.  Assumption~3 in Lemma~\ref{lIR}
holds because of the choice of $U$ (see~\eqref{eqphi1}). Thus,
Lemma~\ref{lIR} implies estimate~\eqref{eqNP0uU}.

Combining Lemma~\ref{lRothe} with estimate~\eqref{eqNP0uU}, we see
that $T_1$ can be chosen equal to $T_0$. Hence, $u\in
W_q^{2,1}(Q_{T_0})$ and, by Theorem~\ref{tLinearHeat}
and~\eqref{eqNP0uU} (with any $T\le T_0$ instead of $T_1$), we
obtain estimates~\eqref{eqNP02*}. Part 1 of the theorem is proved.

2. Let us prove~\eqref{equ-un}. Assume the contrary: there is
$\varepsilon>0$ and a subsequence of $u_n$ (which we denote $u_n$
again) such that
\begin{equation}\label{equ-un1}
\|u_n-u\|_{W_q^{2,1}(Q_T)}\ge \varepsilon,\quad n=1,2,\dots.
\end{equation}

2a. Note that~\eqref{equ0-u0n} and Lemma~\ref{lu0b0v0} imply
\begin{equation}\label{equ-un1'}
\|v_{0n}-v_0\|_{L_q(Q_T)}\to0,\quad n\to\infty.
\end{equation}

Further, by part 1 of the theorem, $u_n$ are bounded in
$W_q^{2,1}(Q_T)$ uniformly with respect to $n$. Therefore, by the
compactness of the embedding $W_q^{2,1}(Q_T)\subset L_q(Q_T)$,
there is a subsequence of $u_n$ (which we denote $u_n$ again) that
is fundamental in $L_q(Q_T)$.

For this subsequence, using Theorem~\ref{tLinearHeat} and the
Lipschitz continuity of $f$ (Condition~\ref{condRegularity}), we
have
\begin{equation}\label{equ-un2}
\begin{aligned}
&\|u_n-u_k\|_{W_q^{2,1}(Q_T)}\\
&\le
k_1\left(\|\phi_n-\phi_k\|_{W_{q,N}^{2-2/q}}+\left(\,\int\limits_{Q_T}|f(u_n,v_{0n})-f(u_k,v_{0k})|^q\,dxdt\right)^{1/q}\right)
\\
&\le k_2\left(\|\phi_n-\phi_k\|_{W_{q,N}^{2-2/q}}+
\left(\,\int\limits_{Q_T}(|u_n-u_k|+|v_{0n}-v_{0k}|)^q\,dxdt\right)^{1/q}\right)\\
&\le
k_2(\|\phi_n-\phi_k\|_{W_{q,N}^{2-2/q}}+\|u_n-u_k\|_{L_q(Q_T)}+\|v_{0n}-v_{0k}\|_{L_q(Q_T)}),
\end{aligned}
\end{equation}
where $k_1,k_2>0$ do not depend on $n$ and $k$. The latter
inequality, the first convergence in~\eqref{equ0-u0n},
relation~\eqref{equ-un1'}, and the fact that $u_n$ is fundamental
in $L_q(Q_T)$ imply that $u_n$ is fundamental and, thus, converges
to some $\hat u$ in $W_q^{2,1}(Q_T)$.

2b. Passing to the limit as $n\to\infty$ and using
Lemma~\ref{lu0b0v0} and Condition~\ref{condRegularity}, we
conclude that $\hat u$ is a solution of the problem
$$
\left\{
\begin{aligned}
& \hat u_t=\hat u_{xx}+f(\hat u,v_0),\quad (x,t)\in Q_T,\\
& \hat u_x|_{x=0}=\hat u_x|_{x=1}=0,\\
& \hat u|_{t=0}=\phi(x),\quad x\in(0,1).
\end{aligned}
\right.
$$
But the latter problem has a unique solution due to part 1 of the
proof. Hence, $\hat u=u$, which contradicts~\eqref{equ-un1}.

Part 2 is proved.

3. Denote
\begin{equation}\label{eqsetOmega}
\Omega=\left\{(x,\phi(x)): x\in[\ob,\ob+1/m], \; \varphi(x) \in
[\alpha,\alpha+1/m^2]\right\}.
\end{equation}

Consider the two cases: $\Omega\ne\varnothing$ and
$\Omega=\varnothing$.

{\bf Case I.} Let $\Omega\ne\varnothing$.
The second inequality in \eqref{eqNP02*} implies that
\begin{equation}\label{eqAdd19.1}
\|u(\cdot, t_1) - u(\cdot, t_2)\|_C + \|u_x(\cdot, t_1) - u_x(\cdot, t_2)\|_C \leq c|t_1-t_2|^{\gamma},
\end{equation}
where $c = c_2(m+f_U)$. This inequality implies that,  for small
enough $\tau(m) > 0$, the functions $\phi$, $\psi_1  = u(\cdot,
t_1)$, $\psi_2  = u(\cdot, t_2)$ with $t_1, t_2 \in [0, \tau(m)]$
satisfy the assumptions of Lemma \ref{lEmA}. Therefore, $a \in
C^{\gamma}$ and
\begin{equation}\label{eqEN4}
\|a\|_{C^{\gamma}} \leq a^* = 1 + 2mc.
\end{equation}

Denote $b(t)=\max\limits_{s\in[0,t]}a(s)$. Note that $b(0) = a(0) = \ob$. It follows from~\eqref{eqAdd19.1} and~\eqref{eqEN4} (decreasing $\tau(m)$ if necessarily) that
\begin{equation}\label{eqEN3'}
b(t)\in[\ob,\ob+1/(2m)], \quad \|b\|_{C^\gamma[0,T]}\le a^*,
\end{equation}
\begin{equation}\label{eqEN6}
u_x(x,t)\ge\dfrac{1}{m+1}\quad \forall
x\in[b(t),b(t)+1/(2m)].
\end{equation}

Denote
$
N=N(m,U)=\max\left(2m, [c_1(m+f_U)]+1\right),
$
where $c_1$ is the constant in~\eqref{eqNP02*} and $[\cdot]$ stands for the integer part of a number.
 Then~\eqref{eqEN3'} and~\eqref{eqEN6} imply
\begin{equation}\label{eqEN7}
a(t),b(t)\in[\ob,\ob+1/N],
\end{equation}
\begin{equation}\label{eqEN8}
u_x(x,t)\ge\dfrac{1}{N}\quad \forall x\in[b(t),b(t)+1/N].
\end{equation}

Now we introduce the function $\xi(x,t)$ as follows: $\xi(x,t)=1$
for $x\le b(t)$ and $\xi(x,t)=2$ for $x> b(t)$. We will see below
that $\xi(x,t)$ is the configuration function of the hysteresis
$H(\xi_0,u)$.

Let us show that $(u(\cdot,t),\xi(\cdot,t)) \in E_N$ on the
interval $t\in[0,T]$, provided that $T = T(m)$ is small enough.

\begin{enumerate}
\item[i.] $b(t)\le 1-1/N$. This follows from~\eqref{eqEN3'}.

\item[ii.] $u(x,t)\le\beta-1/N^2$ for $x\in[0,\ob]$. This follows
for sufficiently small $T$ from the fact that
$u(x,0)=\phi(x)\le\beta-1/m^2$ for $x\in[0,\ob]$ and from the
H\"older continuity of $u$ (the second inequality
in~\eqref{eqNP02*}).

\item[iii.] $u(x,t)\ge\alpha+1/N^2$ for $x\in[\ob+1/N,1]$. This
follows for sufficiently small $T$ from the fact that
$u(x,0)=\phi(x)\ge\alpha+1/(Nm)$ for $x\in[\ob+1/N,1]$ and
from the H\"older continuity of $u$ (the second inequality
in~\eqref{eqNP02*}).

\item[iv.] If $x\in[b(t),b(t)+1/N]$, then $u(x,t)\ge\alpha$ and
$u_x(x,t)\ge 1/N$. The first inequality holds by construction of
the function    $b(t)$. The second inequality follows
from~\eqref{eqEN6}.

\item[v.] $\|u(\cdot,t)\|_{W_{q,N}^{2-2/q}}\le N$. This follows
for sufficiently small $T$ from the fact that
$\|\phi\|_{W_{q,N}^{2-2/q}}\le m$ and from the first inequality
in~\eqref{eqNP02*}.
\end{enumerate}

Items i--v together with relations~\eqref{eqEN4}, \eqref{eqEN3'},
and~\eqref{eqEN7} prove assertions 3.(a)-3.(c) of the theorem in
Case I.

\medskip

{\bf Case II.} Let $\Omega=\varnothing$, where $\Omega$ is the
set in \eqref{eqsetOmega}. Since $u(x,0)=\phi(x) \geq 1/m^2$ for $x\in [\ob,\ob+1/m]$, it follows from   the H\"older
continuity of $u$ (the second inequality in~\eqref{eqNP02*}) that
$u(x,t)\ge 1/(m+1)^2$ for $x\in[\ob,\ob+1/(m+1)]$, provided that
$T$ is sufficiently small.

In this case $a(t)=\ob$ for all $t\in[0,T]$ and it is easy to see
that   parts 1--5 of the definition of $E_N$ again hold for the
pair $(u(\cdot,t),\xi(\cdot,t))$ and for $N=m+1$
on the interval $t\in[0,T]$, provided that $T$ is small enough.

Thus, we have proved assertions 3.(a)-3.(c) of the theorem in Case
II.
\endproof

\begin{remark}\label{rNP0anyu0}
Let
$$
f_0(u,x,t)=f(u,v_0(x,t)),\quad (x,t)\in Q_{T_0},
$$
where $v_0(x,t)$ is the spatially distributed hysteresis defined
for some $u_0\in C(\oQ_T)$:
$$
v_0(x,t)=\cH(\xi_0(x),u_0(x,\cdot))(t).
$$
Let
$$
\|u_0\|_{C(\oQ_{T_0})}\le U.
$$
 and the function $v_0(x,t)$
be measurable. Then parts 1 and 3 of Theorem~\ref{tNP0} remain
true. The proof is analogous to that for Theorem~\ref{tNP0}.
\end{remark}

\section{Proof of the main results}\label{SecProofs}

\subsection{Existence of solutions}

In this subsection, we prove an analogue of
Theorem~\ref{thLocalExistence} under the assumptions made in
Sec.~\ref{subsecSimple}.

As before, we choose $m$ such that $(\phi,\xi_0)\in E_m$. It will
be convenient to reformulate the theorem.

\begin{theorem}[local existence]\label{thLocalExistenceSimple} Let
Condition~$\ref{condA1'}$ be satisfied. Then the following hold
for the number $T=T(m)>0$ from part $3$ of Theorem~$\ref{tNP0}$.
\begin{enumerate}
\item Any solution $u\in W_q^{2,1}(Q_T)$ of
problem~\eqref{eq1}--\eqref{eq3}  in $Q_T$  is transverse and
preserves spatial topology. Moreover, it possesses all the
properties from part $3$ of Theorem~$\ref{tNP0}$.

 \item
There is at least one transverse topology preserving solution
$u\in W_q^{2,1}(Q_T)$ of problem~\eqref{eq1}--\eqref{eq3}
in~$Q_T$.
\end{enumerate}
\end{theorem}

\proof 1. Let $u\in W_q^{2,1}(\oQ_T)$ be an arbitrary solution of
problem~\eqref{eq1}--\eqref{eq3}. Using Lemma~\ref{lIR} (as in the
proof of Theorem~\ref{tNP0}), we obtain that $\|u\|_{C(\oQ_T)}<
U$. Therefore, by Remark~\ref{rNP0anyu0}, the first assertion of
the theorem is true.

2. Let us prove the second assertion.   Let $\cR\subset
C(\oQ_T)\times C[0,T]$ be the set defined before
Lemma~\ref{lu0b0v0}. Clearly, $\cR$ is a closed convex set.

Take any $(u_0,b_0)\in\cR$ and define $f_0(u,x,t)$ by
formula~\eqref{eqf0}. Then Theorem~\ref{tNP0} implies that
problem~\eqref{eqPf0} has a unique transverse topology preserving
solution $u\in W_q^{2,1}(Q_T)\subset C^\gamma(\oQ_T)$.

Consider the functions $a(t)$ and $b(t)$ (see part 3 of
Theorem~\ref{tNP0}). By Theorem~\ref{tNP0}, $a,b\in C^\gamma[0,T]$
and the function $b(t)$ defines the (unique) discontinuity point
of the hysteresis $\cH(\xi_0,u)$ and of its configuration function
$\xi(x,t)$ at each moment $t\in[0,T]$. Moreover $(u,b)\in\cR$.

Thus, we can define a nonlinear operator $R:\cR\to\cR$ by the
formula $R(u_0,b_0)=(u,b)$.

Let us show that $R$ is continuous. Let a sequence
$(u_{0n},b_{0n})$ converge to $(u_0,b_0)$. We define
$f_{0n}(u,x,t)$ by formula~\eqref{eqf0} with
$u_{0n},b_{0n},v_{0n}$ instead of $u_0,b_0,v_0$. Let $u_n$ be a
solution of problem~\eqref{eqPf0} with the right-hand side
$f_{0n}$. By Theorem~\ref{tNP0}, there is a number
$N=N(m,U)\in\mathbb N$ such that
\begin{equation}\label{21*}
(u,\xi),(u_n(\cdot,t),\xi_n(\cdot,t))\in E_N\quad \forall
t\in[0,T],\ n=1,2,\dots.
\end{equation}
Let $a_n(t)$ corresponds to $u_n(x,t)$ in the same way as $a(t)$
corresponds to $u(x,t)$. Set $b_n(t)=\max\limits_{s\in[0,t]}
a_n(s)$ for $t\in[0,T]$. Then
$$
R(u_{0n},b_{0n})=(u_n,b_n).
$$

By the construction of $T$ ($\leq \tau(m)$) in the proof of
assertion 3 of Theorem \ref{tNP0}, the functions $\phi$, $\psi_1 =
u(\cdot, t)$, $\psi_2 = u_n(\cdot, t)$ satisfy the assumptions of
Lemma \ref{lEmA} and, hence,
$$
|a_n(t) - a(t)| \leq 2m\|u_n-u\|_{C(\overline{Q}_T)}.
$$
This implies that $a_n$ converges to $a$ in $C[0, T]$. Thus, Lemma~\ref{lab}
implies that $b_n$ converges to $b$ in $C[0,T]$. Therefore, the
operator $R$ is continuous on $\cR$.

The operator $R$ is also compact. Indeed, the map
$\cR\ni(u_0,b_0)\mapsto u\in C(\oQ_T)$ is compact due
to~\eqref{eqNP02*} and the compactness of the embedding
$C^\gamma(\oQ_T)\subset C(\oQ_T)$. The map $\cR\ni(u_0,b_0)\mapsto
b\in C[0,T]$ is compact due to~\eqref{equ2'}, part 1 of
Lemma~\ref{lab}, and   the compactness of the embedding
$C^\gamma[0,T]\subset C[0,T]$.

Now the Schauder fixed-point theorem implies that  the operator
$R$ has a fixed point $(u,b)\in\cR$. Taking into account
Remark~\ref{rH=v0} and \eqref{21*}, we see that $u$ is a transverse topology
preserving solution of problem~\eqref{eq1}--\eqref{eq3}.
\endproof

\subsection{Continuation of solutions}\label{secCont}

In this subsection, we prove an analogue of Theorem~\ref{tCont}
under the assumptions made in Sec.~\ref{subsecSimple}.

\proof[Proof of Theorem~$\ref{tCont}$]

1. We assume that there is $T_0>0$ such that $u(x,t)$ cannot be
continued to $[0,T_0]$ as a transverse solution of
problem~\eqref{eq1}--\eqref{eq3}.

Applying Theorem~\ref{thLocalExistenceSimple} and using part 1 of
Lemma~\eqref{lEm}, we obtain a sequence $m_k\in\mathbb N$
($k=1,2,\dots$) such that $m_{k+1} > m_k$ and a sequence of time
moments  $t_k$ ($k=1,2,\dots$) such that $t_{k+1} > t_k$ with the
following properties.
\begin{enumerate}
\item For each $k$, the solution $u(x,t)$ of
problem~\eqref{eq1}--\eqref{eq3} can be continued as a transverse
solution to the time interval $[0,t_k]$.

\item For each $k$,
\begin{equation}\label{eqCont1}
u(\cdot, t_k) \in E_{m_k} \setminus E_{m_k-1}.
\end{equation}
\end{enumerate}

Denote $ T=\lim\limits_{k\to\infty}t_k. $ By assumption, $T\le
T_0$.

Since $u$ is a solution of problem~\eqref{eq1}--\eqref{eq3} in
$Q_{t_k}$ for all $k$ and
$$
\|u\|_{W_q^{2,1}(Q_{t_k})}  \le c_1(\|\phi\|_{W_{q,N}^{2-2/q}}+f_U)
$$
by Remark~\ref{rNP0anyu0} (with the right-hand side not depending
on $k$), it follows that $u\in W_q^{2,1}(Q_T)$ and $u$ is a
solution of problem~\eqref{eq1}--\eqref{eq3}  in $Q_T$. Since
$u(\cdot,t)$ is a continuous $W_{q,N}^{2-2/q}$-valued function, we
have
\begin{equation}\label{eqCont2}
\|u(·, t_k)-u(\cdot,T)\|_{W_{q,N}^{2-2/q}}\to 0,\quad k\to\infty.
\end{equation}

Denote by $b(t)$  the discontinuity point of the configuration
function $\xi(x,t)$. By construction $b(t)$ is continuous and
nondecreasing on $[0,t_k]$ for all $k$. Therefore, $b(t)$ is
continuous on $[0,T]$. In particular,
\begin{equation}\label{eqCont3}
b(t_k)-b(T)\to 0,\quad k\to\infty.
\end{equation}

It follows from~\eqref{eqCont1}--\eqref{eqCont3} and from part 2
of Lemma~\ref{lEm} that $b(T)=1$ or $u(x,T)$ is not transverse
with respect to $\xi(x,T)$.

2. Now we consider several cases.

\begin{figure}[ht]
\begin{minipage}[t]{0.60\linewidth}
        \centering
        \includegraphics[width=\linewidth]{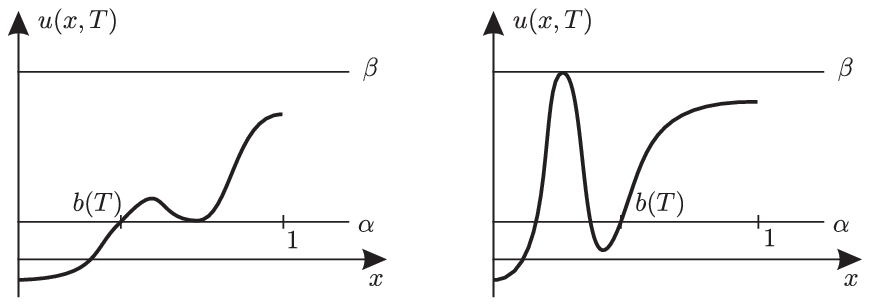}
        \caption{Transversality fails at a point different from $b(T)$}
        \label{fig7}
\end{minipage}
\hspace{0.5cm}
\begin{minipage}[t]{0.29\linewidth}
        \centering
        \includegraphics[width=\linewidth]{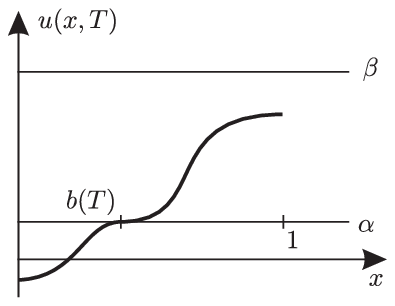}
        \caption{Transversality fails at the point $b(T) \ne 1$}
        \label{fig8}
        \end{minipage}
\end{figure}

2.1. First, we assume that $b(T)<1$. Then $u(x,T)$ is not
transverse with respect to $\xi(x,T)$. This happens because the
graph of the function $u(x,T)$ touches the line $\alpha$ or
$\beta$ (see Fig.~\ref{fig7}), or because $u(b(T),T)=\alpha$,
$u_x(b(T),T)=0$ (see Fig.~\ref{fig8}).

\begin{figure}[ht]
        \centering
        \includegraphics[width = \textwidth]{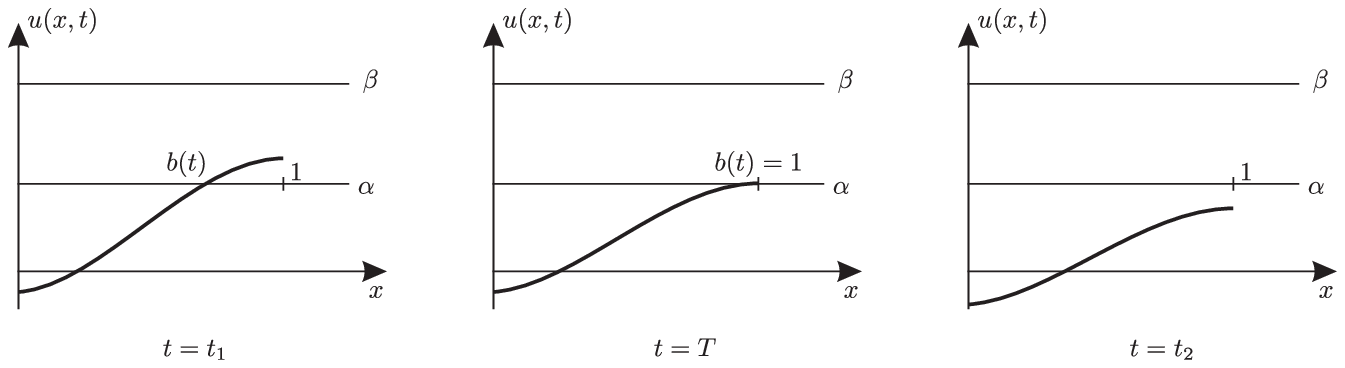}
        \caption{Transversality persists, but topology changes: ${\mathcal H}  = 1$ for $t \ge T$}
        \label{fig9}
\end{figure}

2.2. If  $b(T)=1$, then $u(1,T)=\alpha$. Furthermore, $u_x(1,T)=0$
due to the Neumann boundary condition~\eqref{eq2}. In this case,
either $u(x,T)$ is not transverse (if the transversality fails at
some point $x\in[0,1)$ at the same moment $t=T$, in which case
$T=T_{max}$) or $u(x,T)$ is transverse and $\xi(x,T)\equiv 2$,
$x\in[0,1]$. In the latter case, we can proceed similarly to the
above, but effectively without hysteresis, i.e.,
$$
\cH(\xi_0(x),u(x,\cdot))(t)\equiv H_2(u(x,t)),\quad t\ge T,
$$
(see Fig.~\ref{fig9}).
Thus, one can find $T_{max}$ as in part 1 of the proof.
\endproof

\subsection{Continuous dependence of solutions on initial data}\label{subsecContInitData}

In this subsection we prove an analogue of
Theorem~\ref{tContDepInitData} under the assumptions from
Sec.~\ref{subsecSimple}. The proof will consist of two steps.
First, we shall prove the continuous dependence on small time
intervals and then on the whole interval $[0,T]$.

\subsubsection{Continuous dependence on small time interval}\label{subsubsecSmallInt1}
  Let the initial data $\phi(x)$ and $\xi_0(x)$ satisfy
Condition~\ref{condA1'}. By part 1 of Lemma~\ref{lEm},
$(\phi,\xi_0)\in E_m$ for some $m\in\mathbb N$. Let $\phi_n(x)$
and $\xi_{0n}(x)$ be  other initial data satisfying
Condition~\ref{condA1'}. Suppose that
\begin{equation}\label{eqSmallInt1}
\|\phi_n-\phi\|_{W_{q,N}^{2-2/q}}\to 0, \quad \ob_n-\ob\to 0,\quad
n\to \infty,
\end{equation}
where $\ob_n$ is the discontinuity point of $\xi_{0n}(x)$.

By part 3 of Lemma~\ref{lEm}, there is $n_1=n_1(m)>0$ such that
\begin{equation}\label{eqSmallInt2}
(\phi,\xi_0), (\phi_n,\xi_{0n})\in E_{m+1}\quad \forall n\ge
n_1(m).
\end{equation}
In what follows, we assume $n\ge n_1(m)$. By
Theorem~\ref{thLocalExistenceSimple}, there is $T_1=T_1(m+1)>0$ for
which problem~\eqref{eq1}--\eqref{eq3} has transverse topology
preserving solutions $u,u_n\in W_q^{2,1}(Q_{T_1})$ with  the
initial data $\phi,\xi_0$ and $\phi_n,\xi_{0n}$, respectively.
Moreover, any solution of problem~\eqref{eq1}--\eqref{eq3} in
$Q_{T_1}$ is transverse and preserves topology.

We introduce the functions $a(t)$ and $a_n(t)$ corresponding to
$u$ and $u_n$ as described in part 3 of Theorem~\ref{tNP0}. Then
the discontinuity points of the corresponding configuration
functions $\xi(x,t)$, $\xi_n(x,t)$ are given by
\begin{equation}\label{eq1.1-5}
b(t)=\max\limits_{s\in[0,t]}a(s),\quad
b_n(t)=\max\limits_{s\in[0,t]}a_n(s),\quad t\in[0,T_1].
\end{equation}

\begin{lemma}\label{lContDepInitData1}
Under the above assumptions and the additional assumption that~$u$
is a unique solution of problem~\eqref{eq1}--\eqref{eq3} in
$Q_{T_1}$ with the initial data $\phi,\xi_0$, we have
$$
\|u_n-u\|_{W_q^{2,1}(Q_{T_1})}\to 0,\quad \|b_n-b\|_{C[0,T_1]}\to
0,\quad n\to\infty.
$$
\end{lemma}
\proof 1. Assume to the contrary that there is $\varepsilon>0$
such that
\begin{equation}\label{eq1.1-4'}
\|u_n-u\|_{W_q^{2,1}(Q_{T_1})}\ge\varepsilon,\quad n=1,2,\dots,
\end{equation}
for some subsequence of $u_n$, which we denote $u_n$ again.

Theorem~\ref{thLocalExistenceSimple} implies that, for all
sufficiently large $n$, the functions $u_n$ and $a_n$  are
uniformly bounded in $W_q^{2,1}(Q_{T_1})$ and $C^\gamma[0,T_1]$,
respectively. Therefore, we can choose subsequences of $u_n$ and
$a_n$ (which we denote $u_n$ and $a_n$ again)  such that
\begin{equation}\label{eq1.1-8}
\|u_n-\hat u\|_{C^\gamma(\oQ_{T_1})}\to 0,\ \|(u_n)_x-\hat
u_x\|_{C^\gamma(\oQ_{T_1})}\to 0,\quad n\to\infty,
\end{equation}
\begin{equation}\label{eq1.1-9}
\|a_n-\hat a\|_{C[0,T_1]}\to 0,\quad n\to\infty
\end{equation}
for some function $\hat u\in C^\gamma(\oQ_{T_1})$ with  $\hat
u_x\in C^\gamma(\oQ_{T_1})$ and some function $\hat a\in
C[0,T_1]$.

Denote
\begin{equation*}
\hat b(t)=\max\limits_{s\in[0,t]}\hat a(s), \quad t\in[0,T_1].
\end{equation*}

It follows from~\eqref{eq1.1-8}, \eqref{eq1.1-9}, and Lemma
\ref{lEmA} that  the functions $\hat u$ and $\hat a$ are such that
\begin{equation*}
\cH(\xi(x),\hat u(x,\cdot))(t)=\begin{cases} H_1(\hat u(x,t)), & 0\le x\le  \hat b(t),\\
H_2(\hat u(x,t)), & \hat b(t)<x\le 1.
\end{cases}
\end{equation*}

 Combining~\eqref{eq1.1-9} with Lemma~\ref{lab}, we obtain
\begin{equation}\label{eq1.1-11}
\|b_n-\hat b\|_{C[0,T_1]}\to 0,\quad  n\to\infty.
\end{equation}

2. Now using~\eqref{eqSmallInt1}, \eqref{eq1.1-8},
\eqref{eq1.1-11}, and part 2 of Theorem~\ref{tNP0}, we see that
the sequence of functions $u_n$ satisfying the problems
$$
\left\{
\begin{aligned}
& (u_n)_t=(u_n)_{xx}+f(u_n,\cH(\xi_n,u_n)),\quad (x,t)\in Q_{T_1},\\
& (u_n)_x|_{x=0}=(u_n)_x|_{x=1}=0,\\
& u_n|_{t=0}=\phi_n(x),\quad x\in(0,1),
\end{aligned}
\right.
$$
converges in $W_q^{2,1}(Q_{T_1})$ to $\hat u$ satisfying the
problem
$$
\left\{
\begin{aligned}
& \hat u_t=\hat u_{xx}+f(\hat u,\cH(\xi,\hat u)),\quad (x,t)\in Q_{T_1},\\
& \hat u_x|_{x=0}=\hat u_x|_{x=1}=0,\\
& \hat u|_{t=0}=\phi(x),\quad x\in(0,1).
\end{aligned}
\right.
$$
But the latter problem has a unique solution by assumption. Hence,
$\hat u=u$, which contradicts~\eqref{eq1.1-4'}.

Thus, the first convergence in the lemma is proved. The second
convergence follows from~\eqref{eq1.1-11}.
\endproof

\subsubsection{Continuous dependence on the whole interval without change of
topology}\label{subsubsecContDepWholeIntNoChangeTop}

Now we shall prove Theorem~\ref{tContDepInitData} under the
assumptions from Sec.~\ref{subsecSimple}. In particular, we assume
that $\phi$ and $\xi_0$ satisfy Condition~\ref{condA1'}. By
assumption, problem~\eqref{eq1}--\eqref{eq3} has a unique
transverse topology preserving solution $u$ in $Q_s$ for any $s\le
T$. We denote by $b(t)$ the (unique) discontinuity point of the
corresponding configuration function $\xi(x,t)$.

Further, we assume that $\phi_n, \xi_{0n}$, $n=1,2,\dots$, is a
sequence of other initial data satisfying Condition~\ref{condA1'}
such that
\begin{equation}\label{Eq3.1-1}
\|\phi_n-\phi\|_{W_{q,N}^{2-2/q}}\to 0,\quad \ob_n-\ob\to 0,\quad
n\to\infty,
\end{equation}
where $\ob_n$ is the discontinuity point of $\xi_{0n}$.

We have to show that, for all sufficiently large $n$,
problem~\eqref{eq1}--\eqref{eq3} with the initial data $\phi_n,
\xi_{0n}$ has at least one transverse topology preserving solution
$u_n\in W_q^{2,1}(Q_T)$ and each sequence of such solutions
satisfies
\begin{equation}\label{Eq3.1-3}
\|u_n-u\|_{W_q^{2,1}(Q_T)}\to 0,\quad \|b_n-b\|_{C[0,T]}\to
0,\quad n\to\infty,
\end{equation}
where $b_n(t)$ is the discontinuity point of $\xi_n(x,t)$.

\proof 1. By Lemma~\ref{lEm},  there is $m\in\mathbb N$ such that
$(u(\cdot,t),\xi(\cdot,t))\in E_m$ for all $t\in[0,T]$. Let us fix
such a number $m$. Suppose that~\eqref{Eq3.1-3} does not hold. Let
$\tau$ be the infimum of the set of all  $s\in[0,T]$ such that at
least one of the convergences
\begin{equation}\label{Eq3.1-3'}
\|u_n-u\|_{W_q^{2,1}(Q_s)}\to 0,\quad \|b_n-b\|_{C[0,s]}\to
0,\quad n\to\infty,
\end{equation}
does not hold. By assumption $\tau<T$. On the other hand,
Lemma~\ref{lContDepInitData1} implies that $\tau\ge T_1$. In
particular, this means that
$$
\|u_n(\cdot,\tau-T_1/2)-u(\cdot,\tau-T_1/2)\|_{W_{q,N}^{2-2/q}}\to
0,\quad b_n(\tau-T_1/2)-b(\tau-T_1/2)\to 0,\quad n\to\infty.
$$
Applying Lemma~\ref{lContDepInitData1} again, we obtain the
convergence in~\eqref{Eq3.1-3'} for $s=\min(\tau+T_1/2,T)>\tau$.
This contradiction proves~\eqref{Eq3.1-3}.
\endproof

\subsubsection{Continuous dependence on the whole interval with change of topology}

It remains to prove Corollary~\ref{corContDepInitData} under the
assumptions from Sec.~\ref{subsecSimple}. In particular, we shall
keep notations from the previous
Sec.~\ref{subsubsecContDepWholeIntNoChangeTop}. However, we
additionally assume that $t_{max}>0$ is the number where the
topology of $u$ changes, i.e.,
$$
b(t_{max})=1,\qquad u(1,t_{max})=0,
$$
and that the solutions $u$ and $u_n$ are transverse on the
interval $[0,T]$ with $T\ge t_{max}$.

\begin{lemma}\label{lContDepInitDatatmax}
$ \|u-u_n\|_{W_q^{2,1}(Q_{t_{max}})}\to 0$ and $
\|b-b_n\|_{C[0,t_{max}]}\to 0$ as $n\to\infty$.
\end{lemma}
\proof 1. Let us prove the first convergence in the lemma. Assume
to the contrary that there is $\varepsilon>0$ such that
\begin{equation}\label{eq4-3}
\|u-u_n\|_{W_q^{2,1}(Q_{t_{max}})}\ge\varepsilon,\quad
n=1,2,\dots,
\end{equation}
for some subsequence of $u_n$, which we denote $u_n$ again.

Theorem~\ref{tContDepInitData} implies that
\begin{equation}\label{Eq4-0}
\|u-u_n\|_{W_q^{2,1}(Q_{t_{max}-\delta})}\to 0,\quad
\|b-b_n\|_{C[0,t_{max}-\delta]}\to 0,\quad n\to\infty,
\end{equation}
for all (small) $\delta>0$. Therefore,
\begin{equation}\label{Eq4-2-1}
u_n(x,t)\to u(x,t)\quad \text{a.e. } (x,t)\in Q_{t_{max}},
\end{equation}
\begin{equation}\label{Eq4-2-2}
b_n(t)\to b(t)\quad \text{a.e. } t\in(0,t_{max}),\quad n\to\infty.
\end{equation}

Furthermore, Remark~\ref{rNP0anyu0} and Lemma~\ref{lIR} imply that
the functions $u_n$   are uniformly bounded in
$W_q^{2,1}(Q_{t_{max}})$.  Therefore, we can choose a subsequence
of $u_n$   (which we denote $u_n$ again) converging in
$C(\oQ_{t_{max}})$. Taking into account~\eqref{Eq4-2-1}, we see
that it converges to $u$:
\begin{equation}\label{Eq4-4}
\|u-u_n\|_{C(\oQ_{t_{max}})}\to 0,\quad n\to\infty.
\end{equation}

On the other hand, it follows from~\eqref{Eq4-2-2}, from the
uniform boundedness of $b_n(t)$, and from Lebesgue's dominated
convergence theorem that
\begin{equation}\label{Eq4-5}
\|b-b_n\|_{L_1(0,t_{max})}\to 0,\quad n\to\infty.
\end{equation}
Combining the convergence of initial data, relations~\eqref{Eq4-4}
and~\eqref{Eq4-5}, and part~2 of Theorem~\ref{tNP0}, we conclude
that
\begin{equation}\label{eq4-6}
\|u-u_n\|_{W_q^{2,1}(Q_{t_{max}})}\to 0,\quad n\to \infty,
\end{equation}
for the chosen subsequence. This contradicts~\eqref{eq4-3}.
Therefore, \eqref{eq4-6} holds for the whole sequence $u_n$.

2. Now we prove the second convergence in the lemma. Suppose it
does not hold. Then, due to the second convergence
in~\eqref{Eq4-0},     there is a subsequence of $t_n$ (which we
denote $t_n$ again) and a number $\varepsilon>0$ such that $t_n\to
t_{max}$ and
$$
|b_n(t_n)-b(t_n)|\ge \varepsilon.
$$
Combining this with the fact that $b(t_n)\to 1$, we see that there
is a number $b^*<1$ such that
\begin{equation}\label{eq4-7'}
b_n(t_n)\le b^*.
\end{equation}

On the other hand, since $b(t_{max})=1$, there is a moment
$t^*<t_{max}$ such that
\begin{equation}\label{eq4-7}
b(t^*)>b^*.
\end{equation}
Now, using the monotonicity of $b_n$, the convergence
in~\eqref{Eq4-0} and inequality~\eqref{eq4-7}, we have
$$
b_n(t_n)\ge b_n(t^*)>b^*
$$
for all sufficiently large $n$. This contradicts~\eqref{eq4-7'}.
\endproof

\proof[Proof of Corollary~$\ref{corContDepInitData}$]
Due to Lemma~\ref{lContDepInitDatatmax}, it remains to show that
\begin{equation}\label{Eq4-8}
\|u_n-u\|_{W_q^{2,1}((0,1)\times(t_{max},T))}\to 0,\quad n\to\infty,
\end{equation}
\begin{equation}\label{Eq4-8'}
\|b_n-b\|_{C[t_{max},T]}\to 0,\quad n\to\infty.
\end{equation}

1. Let us prove~\eqref{Eq4-8}. Assume to the contrary that there
is $\varepsilon>0$ such that
\begin{equation}\label{eq4-9}
\|u_n-u\|_{W_q^{2,1}((0,1)\times(t_{max},T))}\ge\varepsilon,\quad
n=1,2,\dots,
\end{equation}
for some subsequence of $u_n$, which we denote $u_n$ again.

Remark~\ref{rNP0anyu0} and Lemma~\ref{lIR} imply that the
functions $u_n$   are uniformly bounded in
$W_q^{2,1}(Q_{t_{max}})$.  Therefore, there is a subsequence of
$u_n$ (which we denote $u_n$ again)  and a function $\hat u$ such
that
\begin{equation}\label{Eq4-10}
\|u_n-\hat u\|_{C([0,1]\times[t_{max},T])}\to 0,\quad n\to\infty.
\end{equation}

Further, we have on the time interval $[t_{max},T]$:
$$
\cH(\xi_{0n},u_n)=\begin{cases}
H_1(u_n), & 0\le x\le  b_n(t),\\
H_2(u_n), & b_n(t)<x\le 1.
\end{cases}
$$
due to the transversality of $u_n$ (here we set $b_n(t)=1$ for
$t>t_{n,max}$ if $b_n(t_{n,max})=1$). Therefore, similarly
to~Lemma~\ref{lu0b0v0}, we obtain
$$
\|\cH(\xi_{0n},u_n)- H_1(\hat u)\|_{L_p((0,1)\times(t_{max},T))}
\le c_0 (T-t_{max})^{1/p}(\|\hat
u-u_n\|_{C([0,1]\times[t_{max},T])}^\sigma+|b_n(t_{max})-1|).
$$
Together with Lemma~\ref{lContDepInitDatatmax} and
relation~\eqref{Eq4-10}, this yields
\begin{equation}\label{Eq4-11}
\|\cH(\xi_{0n},u_n)- H_1(\hat u))\|_{L_p((0,1)\times(t_{max},T))}
\to 0,\quad n\to\infty.
\end{equation}

On the other hand, by Lemma~\ref{lContDepInitDatatmax},
\begin{equation}\label{Eq4-12}
\|u_n|_{t=t_{max}}-u|_{t=t_{max}}\|_{W_{q,N}^{2-2/q}}\to 0,\quad
n\to\infty.
\end{equation}

Using~\eqref{Eq4-10}--\eqref{Eq4-12} and
Theorem~\ref{tLinearHeat}, we see that (cf. proof of part~2 of
Theorem~\ref{tNP0})
\begin{equation}\label{Eq4-8''}
\|u_n-\hat u\|_{W_q^{2,1}((0,1)\times(t_{max},T))}\to 0,\quad
n\to\infty,
\end{equation}
and $\hat u$ is a solution of the problem
$$
\left\{
\begin{aligned}
& \hat u_t=\hat u_{xx}+f(\hat u,H_1(\hat u)),\quad (x,t)\in (0,1)\times(t_{max},T)),\\
& \hat u_x|_{x=0}=\hat u_x|_{x=1}=0,\\
& \hat u|_{t=t_{max}}=u|_{t=t_{max}},\quad x\in(0,1).
\end{aligned}
\right.
$$
Due to the uniqueness assumption, $\hat u=u$ in
$(0,1)\times(t_{max},T))$. Together with~\eqref{Eq4-8''} this
contradicts~\eqref{eq4-9}.

2. To prove~\eqref{Eq4-8'}, we note that, for all
$t\in[t_{max},T]$,
$$
|b_n(t)-b(t)|=b(t_{max})-b_n(t)\le b(t_{max})-b_n(t_{max})\to
0,\quad n\to\infty,
$$
due to~Lemma~\ref{lContDepInitDatatmax}.
\endproof

\section{Some generalizations}\label{secGeneral}

In this section, we  generalize Condition~\ref{condDissip} to the
following one (cf. Remark~\ref{remGeneral} and
Example~\ref{Exh1h2}).

\begin{condition}[generalized dissipativity]\label{condDissipGeneral}
\begin{enumerate}
\item For all sufficiently large $u$,
$$
f(u,H_2(u))\le 0,\quad f(-u,H_1(-u))\ge 0;
$$

\item there is a Lipschitz continuous function $h(u)$ such that
$uh(u)>0$ for $u\ne 0$ and, for any $($small$)$ $\mu>0$, there
exists $U_\mu >0$ such that the function
$$
f_\mu(u,v)=f(u,v)-\mu h(u)
$$
satisfies
$$
f_\mu(U_\mu,H_j(u))< 0,\quad f_\mu(-U_\mu,H_j(u))> 0\quad
\forall |u|\le U_\mu,\ j=1,2.
$$
\end{enumerate}
\end{condition}

Let us prove  Theorem~$\ref{thLocalExistence}$ under
Condition~$\ref{condDissipGeneral}$.

From now on we assume that $(\varphi, \xi_0) \in E_m$ for some $m$ (see Lemma \ref{lEm}) and $0<\mu\le 1$. Due to part 2 of
Condition~\ref{condDissipGeneral}, $f_\mu$ satisfies
Condition~\ref{condDissip}. Therefore, by
Theorems~\ref{thLocalExistence} and~\ref{tCont},
problem~\eqref{eq1}--\eqref{eq3} with the right-hand side $f_\mu$
has a solution $u_\mu$, which can be continued to a maximal
interval of transverse existence $[0,T_{\mu,max})$.

Now the important step is to prove the boundedness of the
solutions $u_\mu$ uniformly with respect to $\mu$.

Let us fix $U>\max(-\alpha,\beta)$ such  that part 1 in
Condition~\ref{condDissipGeneral} holds for all $u\ge U$ and that
$\|\phi\|_{W_{q,N}^{2-2/q}}<U$.

\begin{lemma}\label{lIRGeneral}
The solutions $u_\mu$ satisfy
$\|u_\mu\|_{C(\oQ_{T_{\mu,max}})}<U$.
\end{lemma}
\proof Part 1 of Condition~\ref{condDissipGeneral} and the
assumption $uh(u)>0$ for $u\ne0$ imply the strict inequalities
$$
f_\mu(U,H_2(U))<0,\quad f_\mu(-U,H_1(-U))>0.
$$
Therefore, denoting
$F(x,t)=f_\mu(u_\mu(x,t),\cH(u_\mu(x,\cdot))(t))$, we can proceed
analogously to the proof of Lem\-ma~\ref{lIR} (with obvious
modifications due to another definition of $F$).
\endproof

Let $[0,t_{\mu,max})$ be a maximal interval, on which the solution
$u_\mu$ both remains transverse and preserves spatial topology. We
claim that there is $T>0$ such that $t_{\mu,max}\ge T$ for all
$\mu>0$. Indeed, suppose that  there is a subsequence of
$t_{\mu,max}$ (which we denote $t_{\mu,max}$ again) such that
$t_{\mu,max}\to 0$ as $\mu\to 0$. By Theorem~\ref{tLinearHeat} and
Lemma~\ref{lIRGeneral},
$$
\begin{aligned}
\|u_\mu\|_{W_q^{2,1}(Q_{t_{\mu,max}})}+\max_{t\in[0,t_{\mu,max}]}\|u_\mu(\cdot,t)\|_{W_{q,N}^{2-2/q}}
&\le c_1(\|\phi\|_{W_{q,N}^{2-2/q}}+f_U+h_U),\\
 \|u_\mu\|_{C^\gamma(\overline
Q_{t_{\mu,max}})}+\|(u_\mu)_x\|_{C^\gamma(\overline
Q_{t_{\mu,max}})} &\le c_2 (\|\phi\|_{W_{q,N}^{2-2/q}}+f_U+h_U ),
\end{aligned}
$$
where $f_U$ is defined in~\eqref{eqVfU}, $h_U=\max\limits_{|u|\le
U}|h(u)|$, and   $ c_2>0$   does not depend on
$\mu,t_{\mu,max},\phi$.

The latter estimate shows   that $u_\mu$ remains transverse and
preserves spatial topology for all sufficiently small
$t_{\mu,max}$. This contradiction proves that $t_{\mu,max}\ge
T>0$.

Applying Theorem~\ref{tLinearHeat} and Lemma~\ref{lIRGeneral}
again, we obtain the estimates
\begin{equation}\label{eqNP02**General}
\|u_\mu\|_{C(\oQ_T)}< U,
\end{equation}
\begin{equation}\label{eqNP02*General}
\begin{aligned}
\|u_\mu\|_{W_q^{2,1}(Q_T)}+\max_{t\in[0,T]}\|u_\mu(\cdot,t)\|_{W_{q,N}^{2-2/q}} \le c_1(\|\phi\|_{W_{q,N}^{2-2/q}}+f_U+h_U),\\
\|u_\mu\|_{C^\gamma(\oQ_T)}+\|(u_\mu)_x\|_{C^\gamma(\oQ_T)}  \le
c_2 (\|\phi\|_{W_{q,N}^{2-2/q}}+f_U+h_U),
\end{aligned}
\end{equation}
where $c_1,c_2>0$ do not depend on $\mu$.

For each $u_\mu$, we have a corresponding function $a_\mu(t)$ (cf.
part 3 of Theorem~\ref{tNP0}). It follows
from~\eqref{eqNP02*General} and from the fact that $a_\mu\in
C^\gamma[0,T]$ that one can choose  subsequences,  which we denote
$u_\mu$ and $a_\mu$ again, converging to some functions $u(x,t)$
and $a(t)$ in $C(\oQ)$ and $C[0,T]$, respectively.

Using~\eqref{eqNP02*General}, one can show that the function
$a(t)$ corresponds to $u(x,t)$ in the same sense as the functions
$a_\mu(t)$ correspond to $u_\mu(x,t)$. In particular, this means
that the corresponding hysteresis operators are given by
$$
\cH(\xi_0,u)=\begin{cases}
H_1(u), & 0\le x\le  b(t),\\
H_2(u), & b(t)<x\le 1,
\end{cases}\qquad
\cH(\xi_0,u_\mu)=\begin{cases}
H_1(u_\mu), & 0\le x\le  b_\mu(t),\\
H_2(u_\mu), & b_\mu(t)<x\le 1,
\end{cases}
$$
where $b(t)=\max\limits_{s\in[0,t]}a(s)$ and
$b_\mu(t)=\max\limits_{s\in[0,t]}a_\mu(s)$.

By  Lemma~\ref{lab}, $\|b_\mu-b\|_{C[0,T]}\to 0$ as $\mu\to 0$.
Thus, applying Remark~\ref{rH=v0} and Lemma~\ref{lu0b0v0}, we see
that $\cH(\xi_0,u_\mu)$ converges to $\cH(\xi_0,u)$ in $L_q(Q_T)$.
Now the Lipschitz continuity of $f$ and the estimate
$|h(u_\mu)|\le h_U$ imply that $f_\mu(u_\mu,\cH(\xi_0,u_\mu))$
converges to $f(u,\cH(\xi_0,u))$ in $L_q(Q_T)$.

Denote by $\hat u$ the solution of the linear parabolic problem
$$
\hat u_t=\hat u_{xx}+f(u,\cH(\xi_0,u)),\quad (x,t)\in Q_T,
$$
$$
\hat u_x|_{x=0}=\hat u_x|_{x=1}=0,
$$
$$
\hat u|_{t=0}=\phi(x),\quad x\in(0,1).
$$
By Theorem~\ref{tLinearHeat}, $u_\mu\to\hat u$ in
$W_q^{2,1}(Q_T)$, hence in $C(\oQ_T)$. Therefore, $u=\hat u$, $u$
is a solution of problem~\eqref{eq1}--\eqref{eq3} with the
right-hand side $f$, and estimates~\eqref{eqNP02**General}
and~\eqref{eqNP02*General} yield the same estimates for $u$.
Theorem~\ref{thLocalExistence} is proved.

The number $T$ which we have obtained above depends only on $m$.
Therefore, the continuation theorem (Theorem~\ref{tCont}) and the
theorem on continuous dependence of solutions on initial data
(Theorem~\ref{tContDepInitData}) under
Condition~\ref{condDissipGeneral} are  proved similarly to
Secs.~\ref{secCont} and~\ref{subsecContInitData} (with the help of
estimates~\eqref{eqNP02**General} and~\eqref{eqNP02*General},
which now hold with $u$). The proof of Theorem~\ref{tUniqueness}
does not depend on the dissipativity condition at all
(see~\cite{GurTikhUniq}).

\bigskip

{\bf Acknowledgement:} The authors are grateful to Willi J\"ager
for drawing their attention to the field of hysteresis and to
Bernold Fiedler and Alexander Nazarov for fruitful discussions.
The research of the first author was supported by the DFG project
SFB 910, by the DAAD program G-RISC, and by the RFBR (project
10-01-00395-a). The research of the third author was supported by
the Alexander von Humboldt Foundation.


\begin{thebibliography}{10}

\bibitem{Alt}
H. W. Alt,  {\it On the thermostat problem}.  Control Cyb., {\bf
14}, 171--193 (1985).


\bibitem{Sobolevskii}
A. Ashyralyev, P. E. Sobolevskii. {\it Well-posedness of Parabolic
Difference Equations}. Operator Theory: Advances and Applications,
69. Birkh\"auser Verlag, Basel (1994).

\bibitem{Evans}
L. C. Evans, M. Portilheiro, {\it Irreversibility and hysteresis
for a forward-backward diffusion equation}.  Math. Models Methods
Appl. Sci., {\bf 14}, no. 11, 1599--1620 (2004).

\bibitem{GurTikhUniq}
P. Gurevich, S. Tikhomirov, {\it Uniqueness of transverse solutions 
for reaction-diffusion equations with spatially distributed hysteresis}.
Preprint.

\bibitem{Jaeger1}
F. C. Hoppensteadt, W. J\"ager, {\it Pattern formation by
bacteria} Lecture Notes in Biomathematics {\bf 38}, 68--81 (1980).

\bibitem{Jaeger2}
F. C. Hoppensteadt, W. J\"ager, C. Poppe, {\it A hysteresis model
for bacterial growth patterns}. Modelling of Patterns in Space and
Time, Lecture Notes in Biomath. {\bf 55} (Springer), 123--134
(1984).


\bibitem{Ilin}
A. M. Il'in, B. A. Markov, {\it Nonlinear diffusion equation and
Liesegang rings}. Doklady Akademii Nauk, {\bf 440}, No. 2,
164--167 (2011); English translation: Doklady Mathematics, {\bf
84}, No. 2, 730--733 (2011).


\bibitem{Kopfova}
J. Kopfova,  {\it Hysteresis in biological models}. Journal of
Physics: Conference Series, {\bf 55},  130--134 (2006).







\bibitem{KrasnBook}
M. A. Krasnosel'skii, A. V. Pokrovskii. {\it Systems with
Hysteresis}. Springer-Verlag. Berlin--Heidelberg--New York (1989).
Translated from Russian: {\it Sistemy s Gisterezisom}. Nauka.
Moscow (1983).

\bibitem{LadSolUral} O. A. Ladyzhenskaya, V. A. Solonnikov, N. N.
Uraltseva. {\it Linear and Quasilinear Equations of Parabolic
Type}. Nauka, Moscow, 1967; English translation: Amer. Math. Soc.,
Providence, RI, 1968.

\bibitem{Plotnikov}
P. I. Plotnikov, {\it Passing to the limit with respect to the
viscosity in an equation with variable parabolicity direction},
Differential Equations, {\bf 30}, 614--622 (1994).

\bibitem{Rothe}
F. Rothe.
{\it Global Solutions of Reaction-Diffusion Systems}. Lecture
Notes in Mathematics, 1072. Springer-Verlag, Berlin (1984).

\bibitem{Smoller}
J. Smoller. {\it Shock Waves and Reaction-Diffusion Equations}.
Second edition. Grundlehren der Mathematischen Wissenschaften
[Fundamental Principles of Mathematical Sciences], 258.
Springer-Verlag, New York (1994).



\bibitem{Triebel}
H. Triebel. {\it Interpolation Theory, Function Spaces,
Differential Operators}. Second edition. Johann Ambrosius Barth,
Heidelberg  (1995).

\bibitem{Visintin}
A. Visintin  {\it Differential Models of Hysteresis}.
Springer-Verlag. Berlin --- Heidelberg (1994).

\bibitem{VisintinSpatHyst}
A. Visintin   {\it Evolution problems with hysteresis in the
source term}. SIAM J. Math. Anal, {\bf 17}, 1113--1138 (1986).



\end{thebibliography}
\end{document}